\documentclass[12pt,a4paper]{article}

\usepackage{a4wide}

\usepackage{epstopdf,amsmath,authblk}% To incorporate .eps illustrations using PDFLaTeX, etc.
\usepackage{subfigure}% Support for small, `sub' figures and tables
\usepackage[british]{babel}

\usepackage[numbers]{natbib}% Citation support using natbib.sty
\usepackage{caption}
\allowdisplaybreaks

\usepackage{amsmath,amssymb,calc,dsfont}
\usepackage{graphicx,array}
\usepackage{tikz}
\usepackage{color} %Mag later weg, geeft voor mij aan waar ik nog even naar moet kijken
\usepackage{enumerate}
\usepackage[toc,page]{appendix}
\usepackage[charter]{mathdesign}
\usepackage{listings}
\usepackage{verbatim}

\definecolor{darkgreen}{rgb}{0,0.6,0}

\usepackage{algorithm}
\usepackage{algpseudocode}
\usepackage{tikz}
\usepackage[margin=1 cm]{caption}
\newcolumntype{H}{>{\setbox0=\hbox\bgroup}c<{\egroup}@{}}

\graphicspath{{./},{./images/}}

\newtheorem{assumption}{Assumption}

\newtheorem{remark}{Remark}

\numberwithin{equation}{section}
\providecommand{\href}[2]{#2}

\newcommand{\beq}{\begin{equation}}
	\newcommand{\eq}{\end{equation}}

\newcommand{\car}[3][0]{
	\begin{scope}[shift={#2}, rotate=#1, scale=0.04]
		\draw[#3,fill=#3](69.5,15.5)--(69.5,16.8)--(69.3,19.3)--(68.7,21.8)--(67.2,24.3)--(64.9,26.6)--(61.9,28.1)--(58.4,29.2)--(55.,29.4)--(52.7,29.2)--(50.6,28.8)--(43.1,28.9)--(41.9,30.9)--(40.5,31.)--(39.3,31.)--(39.5,29.)--(28.5,28.7)--(18.4,29.8)--(13.7,29.8)--(10.,29.5)--(2.8,27.9)--(1.7,26.9)--(0.4,23.)--(0.,15.5)--(0.,15.5)--(0.4,8.)--(1.7,4.1)--(2.8,3.1)--(10.,1.5)--(13.7,1.2)--(18.4,1.2)--(28.5,2.3)--(39.5,2.)--(39.3,0.)--(40.5,0.)--(41.9,0.1)--(43.1,2.1)--(50.6,2.2)--(52.7,1.8)--(55.,1.6)--(58.4,1.8)--(61.9,2.9)--(64.9,4.4)--(67.2,6.7)--(68.7,9.2)--(69.3,11.7)--(69.5,14.2)--(69.5,15.5);
		\draw[white,fill=white](53.,15.5)--(52.8,18.1)--(52.,20.4)--(51.,22.2)--(49.8,23.8)--(48.6,24.9)--(47.3,25.1)--(45.7,24.9)--(43.,24.6)--(40.2,23.6)--(39.,22.7)--(38.5,21.1)--(38.7,15.5)--(38.7,15.5)--(38.5,9.9)--(39.,8.3)--(40.2,7.4)--(43.,6.4)--(45.7,6.1)--(47.3,5.9)--(48.6,6.1)--(49.8,7.2)--(51.,8.8)--(52.,10.6)--(52.8,12.9)--(53.,15.5);
		\draw[white,fill=white](17.6,15.5)--(17.7,20.4)--(17.5,21.7)--(16.6,22.9)--(15.3,23.4)--(8.9,22.2)--(7.7,20.9)--(6.9,19.8)--(6.3,18.3)--(5.9,15.5)--(5.9,15.5)--(6.3,12.7)--(6.9,11.2)--(7.7,10.1)--(8.9,8.8)--(15.3,7.6)--(16.6,8.1)--(17.5,9.3)--(17.7,10.6)--(17.6,15.5);
		\draw[white,fill=white](64.5,25.3)--(65.3,24.3)--(66.4,23.2)--(67.1,22.1)--(67.5,20.8)--(66.2,20.8)--(65.2,21.5)--(64.2,22.8)--(63.7,24.2)--(63.5,26.)--(64.5,25.3);
		\draw[white,fill=white](63.5,5.7)--(62.5,5.)--(62.7,6.8)--(63.2,8.2)--(64.2,9.5)--(65.2,10.2)--(66.5,10.2)--(66.1,8.9)--(65.4,7.8)--(64.3,6.7)--(63.5,5.7);
	\end{scope}
}

\title{The Fixed-Cycle Traffic-Light queue with multiple lanes and temporary blockages}

\newcommand{\pedestrian}[3][0]{
	\begin{scope}[shift={#2}, rotate=#1, scale=0.06667]
		\draw[#3,fill=#3](45,16.8)--(44,18)--(43,19)--(42,19.8)--(41,20.3)--(40,20.6)--(39,20.8)--(37.5,21)--(36,20.8)--(35,20.6)--(34,20.3)--(33,19.8)--(32,19)--(31,18)--(30,16.8)--(30,15.5)--(31,14.3)--(32,13.3)--(33,12.5)--(34,12)--(35,11.7)--(36,11.5)--(37.5,11.3)--(39,11.5)--(40,11.7)--(41,12)--(42,12.5)--(43,13.3)--(44,14.3)--(45,15.5);
	\end{scope}
}
\usepackage{tikz}
\usetikzlibrary{decorations.pathreplacing}
\usetikzlibrary{positioning}
\newtheorem{example}{Example}

\author[1]{Rik W. Timmerman}
\author[1]{Marko A. A. Boon}
\affil[1]{Eindhoven University of Technology}

\begin{document}
	%%%%%%%%%%%%%%%%%%%%%%%%%%%%%%%%%%%%%%%%%%
	
	%%%%%%%%%%%%%%%%%%%%%%%%%%%%%%%%%%%%%%%%%%
	%\setcounter{section}{-1} %% Remove this when starting to work on the template.
	%\section{How to Use this Template}
	%The template details the sections that can be used in a manuscript. Note that the order and names of article sections may differ from the requirements of the journal (e.g., the positioning of the Materials and Methods section). Please check the instructions for authors page of the journal to verify the correct order and names. For any questions, please contact the editorial office of the journal or support@mdpi.com. For LaTeX related questions please contact latex@mdpi.com.
	%The order of the section titles is: Introduction, Materials and Methods, Results, Discussion, Conclusions for these journals: aerospace,algorithms,antibodies,antioxidants,atmosphere,axioms,biomedicines,carbon,crystals,designs,diagnostics,environments,fermentation,fluids,forests,fractalfract,informatics,information,inventions,jfmk,jrfm,lubricants,neonatalscreening,neuroglia,particles,pharmaceutics,polymers,processes,technologies,viruses,vision
	\maketitle
	
	\begin{abstract}
		\noindent Traffic-light modelling is a complex task, because many factors have to be taken into account. In particular, capturing all traffic flows in one model can significantly complicate the model. Therefore, several realistic features are typically omitted from most models. We introduce a mechanism to include pedestrians and focus on situations where they may block vehicles that get a green light simultaneously. More specifically, we consider a generalization of the Fixed-Cycle Traffic-Light (FCTL) queue. Our framework allows us to model situations where (part of the) vehicles are blocked, e.g. by pedestrians that block turning traffic and where several vehicles might depart simultaneously, e.g. in case of multiple lanes receiving a green light simultaneously. We rely on probability generating function and complex analysis techniques which are also used to study the regular FCTL queue. We study the effect of several parameters on performance measures such as the mean delay and queue-length distribution.
	\end{abstract}
	
%	\begin{keywords} Traffic-light modelling, Fixed-Cycle Traffic-Light queue, pedestrian crossings, capacity analysis, queueing theory.
%	\end{keywords}
	
	\section{Introduction}
	
	Traffic lights are currently omnipresent in urban areas and one of their aims is to let vehicles drive across an intersection in such a way that the delay is as small as possible. The modelling of queues in front of traffic lights therefore has always been and still is an important topic of study in road-traffic engineering. The overall aim is to create a model that is as realistic as possible, which poses to be a difficult task. There are many studies devoted to traffic control at intersections, ranging from simulation studies and the use of artificial intelligence to analytical and explicit calculations to find good control strategies. This study provides a more realistic extension of the so-called Fixed-Cycle Traffic-Light (FCTL) queue, see e.g.~\cite{darroch1964traffic}, which allows us to perform analytical computations. We call the model that we consider in this paper the blocked Fixed-Cycle Traffic-Light (bFCTL) queue with multiple lanes. Our main aim is to provide an exact computation of the steady-state queue length of the bFCTL queue with multiple lanes, although a transient analysis (possibly with time-varying parameters) is also possible.
	
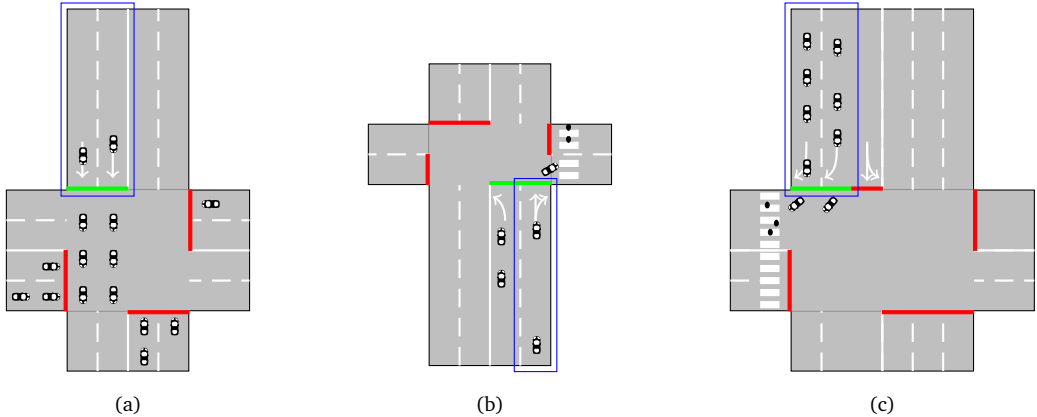
\begin{figure}[h!]
		\centering
		\begin{tabular}{ccc}
			\begin{tikzpicture}[scale=0.08]
				\draw[black,fill=lightgray](30,20) rectangle (50,80);
				
				\draw[black,fill=lightgray](20,30) rectangle (60,50);
				\draw[thick,white,dash pattern=on 7 off 4](20,45) to (60,45);
				\draw[thick,white,dash pattern=on 7 off 4](20,35) to (60,35);
				\draw[thick,white,dash pattern=on 7 off 4](45,20) to (45,80);
				\draw[thick,white,dash pattern=on 7 off 4](35,20) to (35,80);
				\draw[thick,white](20,40) to (60,40);
				\draw[thick,white](40,20) to (40,80);
				\draw[lightgray,fill=lightgray](30,30) rectangle (50,50);
				\draw[red,fill=red](29.5,30) rectangle (30,40);
				\draw[red,fill=red](50.5,50) rectangle (50,40);
				\draw[red,fill=red](50,29.5) rectangle (40,30);
				\draw[green,fill=green](30,50) rectangle (40,50.5);
				\draw[blue](29,49) rectangle (41,81);
				
				\draw[white,thick,->] (32.5,58) to (32.5,52);
				\draw[white,thick,->] (37.5,58) to (37.5,52);
				
				\car{(21,31.7)}{black}
				\car{(26,31.7)}{black}
				\car{(26,36.7)}{black}
				\car[90]{(48.3,26)}{black}
				\car[90]{(43.3,26)}{black}
				\car[90]{(43.3,21)}{black}
				
				\car[270]{(32, 46)}{black}
				\car[270]{(32, 40)}{black}
				\car[270]{(32, 34)}{black}
				\car[270]{(37, 40)}{black}
				\car[270]{(37, 34)}{black}
				\car[270]{(37, 46)}{black}
				\car[270]{(37, 59)}{black}
				\car[270]{(32, 57)}{black}
				
				\car[180]{(55, 48.3)}{black}
			\end{tikzpicture} & \hspace{1 cm}
			\begin{tikzpicture}[scale=0.08,rotate=90]
				\draw[black,fill=lightgray](30,20) rectangle (40,60);
				
				\draw[black,fill=lightgray](0,30) rectangle (50,50);
				\draw[thick,white,dash pattern=on 7 off 4](0,35) to (50,35);
				\draw[thick,white](0,40) to (50,40);
				\draw[thick,white,dash pattern=on 7 off 4](0,45) to (50,45);
				\draw[thick,white,dash pattern=on 7 off 4](35,20) to (35,60);
				\draw[lightgray,fill=lightgray](30,30) rectangle (40,50);
				\draw[white,->,thick](24,32.5) to [out = 0, in = 120] (28,30.5);
				\draw[white,->,thick](24,32.5) to [out = 0, in = 180] (28.5,32.5);
				\draw[white,->,thick](24,37.5) to [out = 0, in = 240] (28,39.5);
				\draw[green,fill=green](30,30) rectangle (30.5,35);
				\draw[green,fill=green](30,35) rectangle (30.5,40);
				\draw[red,fill=red](40,40) rectangle (40.5,50);
				\draw[red,fill=red](30,50) rectangle (35,50.5);
				\draw[red,fill=red](35,30) rectangle (40,30.5);
				\car{(2,31.7)}{black}
				\car{(21,31.7)}{black}
				\car{(20,37.5)}{black}
				\car{(13,37.5)}{black}
				\draw[white,fill=white](31,25.5) rectangle (32,28.5);
				\draw[white,fill=white](33,25.5) rectangle (34,28.5);
				\draw[white,fill=white](36,25.5) rectangle (37,28.5);
				\draw[white,fill=white](38,25.5) rectangle (39,28.5);
				\pedestrian[0]{(37,26)}{black}
				\pedestrian[0]{(35,26)}{black}
				\car[300]{(31.3,31)}{black}	
				\draw[blue](-1,29) rectangle (31,36);
			\end{tikzpicture} & \hspace{1 cm}
			\begin{tikzpicture}[scale=0.08]
				\draw[black,fill=lightgray](30,20) rectangle (60,80);
				\draw[black,fill=lightgray](20,30) rectangle (70,50);
				\draw[thick,white,dash pattern=on 7 off 4](20,40) to (70,40);
				\draw[thick,white,dash pattern=on 7 off 4](20,35) to (70,35);
				\draw[thick,white,dash pattern=on 7 off 4](45,20) to (45,90);
				\draw[thick,white,dash pattern=on 7 off 4](35,20) to (35,90);
				\draw[thick,white,dash pattern=on 7 off 4](50,20) to (50,90);
				\draw[thick,white,dash pattern=on 7 off 4](55,20) to (55,90);
				\draw[thick,white](20,40) to (70,40);
				\draw[thick,white](45,20) to (45,80);
				\draw[lightgray,fill=lightgray](30,30) rectangle (60,50);
				\draw[red,fill=red](29.5,30) rectangle (30,40);
				\draw[red,fill=red](60.5,50) rectangle (60,40);
				\draw[red,fill=red](60,29.5) rectangle (45,30);
				\draw[green,fill=green](30,50) rectangle (40,50.5);
				\draw[red,fill=red](40,50) rectangle (45,50.5);
				
				\draw[blue](29,49) rectangle (41,81);
				
				\draw[white,->,thick](32.5,58) to [out = 270, in = 30] (30.5,52);
				\draw[white,->,thick](37.5,58) to [out = 270, in = 30] (35.5,52);
				\draw[white,->,thick](42.5,58) to (42.5,51.5);
				\draw[white,->,thick](42.5,58) to [out=270, in = 150] (44.5,52);
				
				\draw[white,fill=white](25,48.5) rectangle (28,49.5);
				\draw[white,fill=white](25,46.5) rectangle (28,47.5);
				\draw[white,fill=white](25,44.5) rectangle (28,45.5);
				\draw[white,fill=white](25,42.5) rectangle (28,43.5);
				\draw[white,fill=white](25,40.5) rectangle (28,41.5);
				\draw[white,fill=white](25,39.5) rectangle (28,38.5);
				\draw[white,fill=white](25,37.5) rectangle (28,36.5);
				\draw[white,fill=white](25,35.5) rectangle (28,34.5);
				\draw[white,fill=white](25,33.5) rectangle (28,32.5);
				\draw[white,fill=white](25,30.5) rectangle (28,31.5);
				
				\car[270]{(32, 76)}{black}
				\car[270]{(32, 70)}{black}
				\car[270]{(32, 64)}{black}
				\car[270]{(37, 75)}{black}
				\car[270]{(37, 66)}{black}
				\car[270]{(37, 60)}{black}
				\car[230]{(36.8, 49)}{black}
				\car[220]{(31.5, 49)}{black}
				\car[270]{(32, 55)}{black}
				\pedestrian[270]{(25,50)}{black}
				\pedestrian[270]{(26.5,47)}{black}
				\pedestrian[270]{(25.5,45.5)}{black}
			\end{tikzpicture}
			\\
			\scriptsize (a)  %Platoon forming approaching the intersection
			& \hspace{1cm}
			\scriptsize (b)  %Crossing the intersection in platoons
			& \hspace{1cm}
			\scriptsize (c)  %Crossing the intersection in platoons
		\end{tabular}
		\caption{A visualization of three intersections that can be modeled by the bFCTL queue with multiple lanes. In~(a), the blue rectangle indicates a combination of lanes which can be analyzed as a bFCTL queue with two lanes. The other lanes at the intersection, the complement of the blue rectangle, can be considered separately because of the fixed settings. In~(b), the blue rectangle indicates a lane that can be modeled as a bFCTL queue with a single lane with blockages. In~(c), the blue rectangle indicates two lanes that we can model as a bFCTL queue with \emph{two} lanes where vehicles are potentially blocked by pedestrians.}
		\label{fig:vis}	
	\end{figure}

	The regular FCTL queue is a well-studied model in traffic engineering, see~\cite{boon2019pollaczek,boon2018networks,darroch1964traffic,hagen1989comparison,mcneil1968solution,newell1965approximation,oblakova2019exact,van2006delay,webster1958traffic}. The typical features of the FCTL queue are:
\begin{itemize}
\item A fixed cycle length, fixed green and red times;
\item A general arrival process;
\item Constant interdeparture times of queued vehicles;
\item Whenever the queue becomes empty during a green period, it remains empty since newly arriving vehicles pass the crossing at full speed without experiencing any delay.
\end{itemize}
Due to all the fixed settings, the model focuses on a single lane and does not capture any dependencies or interactions with other lanes.
Unfortunately, in many cases the FCTL queue cannot be applied as a realistic model to study the queue-length distribution in front of a traffic light. Take, for example, an intersection where vehicles from a single stream are spread onto two lanes which are both heading straight and where both lanes are governed by the same traffic light, see also Figure~\ref{fig:vis}(a). Indeed, since there are two parallel lanes in each direction, two vehicles can cross the intersection simultaneously and vehicles will in general switch lanes (if needed) to join the lane with the shorter queue. Moreover, it might be the case that the vehicles are blocked during the green period, e.g. because of a pedestrian crossing the intersection (receiving a green light at the same time as the stream of vehicles that we model), see Figure~\ref{fig:vis}(b) for a visualization. Such blockages might also occur in a multi-lane scenario (where all lanes are going in the same direction) as visualized in Figure~\ref{fig:vis}(c). It is apparent that these situations cannot be modeled by the standard FCTL queue. However, it is extremely relevant to understand such intersections better as is also indicated in e.g.~\cite{tageldin2019models,yan2018design} and more generally, it is e.g. important to investigate pedestrian behaviour at intersections as is done in e.g.~\cite{zhou2019simulation}.
	%Traffic lights are currently omnipresent in urban areas and aim at facilitating traffic in such a way that the traffic delay is as small as possible.
	%The modelling of queues in front of traffic lights is and has been an important topic of study in road-traffic engineering as we have seen before. The overall aim is to create a model that is as realistic as possible, which poses to be a difficult task. %There are many studies devoted to traffic control at intersections, ranging from simulation studies and the use of artificial intelligence to analytical and explicit calculations to find good control strategies.
	The study in this paper provides an extension of the FCTL queue to account for such situations. They seem to be the most common in practice, see e.g.~\cite{shaoluen2020random} for another study on the case as in Figure~\ref{fig:vis}(b). For extensions and other scenarios, we refer the reader to Section~\ref{sec:discussion}. Note that the blocking mechanisms discussed in this paper give rise to more complicated model dynamics and dependencies, which make it impossible to use traditional methods (e.g. Webster's approximation for the mean delay \cite{webster1958traffic}).
	
	A shared right-turn lane as in Figure~\ref{fig:vis}(b), that is a lane with vehicles that are either turning right or are heading straight, has been studied before. However, to the best of our knowledge, there are no papers with a rigorous analysis taking stochastic effects into account while computing e.g. the mean queue length for such lanes. Shared right-turn lanes where vehicles are blocked by pedestrians crossing immediately after the right turn have been considered in e.g.~\cite{alhajyaseen2013left,chen2011saturation,chen2014investigation,chen2008influence,shaoluen2020random,milazzo1998effect,roshani2017effect,rouphail1997pedestrian}. Several case studies, such as~\cite{chen2014investigation} and \cite{roshani2017effect}, indicate that there is a potentially severe impact by pedestrians blocking vehicles. This is for example also reflected in the Highway Capacity Manual (HCM) as published by the Transportation Research Board~\cite{manual2010}, where the focus is on capacity estimation. Most papers have also focused on the estimation of the so-called saturation flow rate, or capacity, of shared lanes where turning vehicles are possibly blocked by pedestrians, see e.g.~\cite{chen2008influence,milazzo1998effect,rouphail1997pedestrian}. In~\cite{chen2011saturation}, it is stated that the used functions for the capacity estimation for turning lanes (such as those in the HCM) might have to be extended to account for stochastic behaviour. In a small case study, \cite{chen2011saturation} confirm that the capacity estimation by the HCM yields an overestimation in various cases. The overestimation of the capacity by the HCM is also observed in several other papers, such as in~\cite{chen2014investigation,chen2008influence} and \cite{shaoluen2020random}, and is probably due to random/stochastic effects. The bFCTL queue explicitly models such stochastic behaviour.
	
	%Pedestrians are not the only potential cause for vehicles being blocked. Other road users might also block (turning) vehicles. The impact of such blockages by other road users has been studied as well. A part of the literature is devoted to the influence of bicycles blocking vehicles, see e.g.~\cite{allen1998effect,chen2018evaluating,chen2007influence,guo2012effect} and references therein. Another prominent line of research is constituted by vehicles that are blocked by vehicles from an opposing stream (for example in case of a non-protected left turn at a signalized intersection). Examples of such studies can be found in e.g.~\cite{chai2014traffic,levinson1989capacity,liu2011arterial,liu2008lane,ma2017two,wu2011modelling,yang2018analytical,yao2013optimal}. Also the effect of shared and exclusive turning lanes is investigated for various scenarios, see for example~\cite{kikuchi2007lengths,tian2006probabilistic,wu1999capacity,zhang2008modeling}.

	A potential application of the bFCTL queue with a single lane as depicted in Figure~\ref{fig:vis}(b) can be found in the model that is studied in~\cite{shaoluen2020random}, which has also been the source of inspiration for this paper. A description of the model in~\cite{shaoluen2020random} is as follows, where we replace the left-turn assumption for left-driving traffic to a right-turn assumption for the more standard case of right-driving traffic. We have a shared lane with straight-going and right-turning traffic controlled by a traffic light, where immediately after the right turn there is a crossing for pedestrians. The pedestrians may block the right-turning vehicles as the vehicles and pedestrians may receive a green light simultaneously. The right-turning vehicles that are blocked, immediately block all vehicles behind them.
	
	Another potential application of the bFCTL queue is to account for bike lanes. Bikes might make use of a dedicated lane or mix with other traffic and in both cases a turning vehicle might be (temporarily) blocked by bicycles because the bicycles happen to be in between the vehicle and the direction that the vehicle is going. As such, blockages have an influence on the performance measures of the traffic light. It is important to take such influences into account in order to find good traffic-light settings. Several papers studying the impact of bikes can be found in~\cite{allen1998effect,chen2007influence,guo2012effect} and \cite{chen2018evaluating}. Also other types of blocking might occur, such as by a shared-left turn lane and opposing traffic receiving a green light simultaneously, see e.g.~\cite{chai2014traffic,levinson1989capacity,liu2011arterial,liu2008lane,ma2017two,wu2011modelling,yang2018analytical,yao2013optimal}. As such, the bFCTL queue (either with multiple lanes or not) is a relevant addition to the literature because it enables a more suitable modelling of traffic lights at intersections with crossing pedestrians and bikes, which leads to traffic-light control strategies for more realistic situations. In order to model a situation where two opposing streams of vehicles potentially block one another as in e.g.~\cite{yang2018analytical}, the bFCTL queue would have to be extended. For more references on the topics discussed in this paragraph see also the review paper by~\cite{cheng2016review}. Another related study is~\cite{oblakova2019exact} who introduce a model with ``distracted'' drivers, which can be considered as an FCTL queue with independent blockages, but this blocking mechanism is a special case of the one discussed in the present paper.

	As mentioned before, we call the model that we consider in this paper the bFCTL queue with multiple lanes. On the one hand we thus allow for the modelling of vehicle streams that are spread over multiple lanes and on the other hand we allow for vehicles to be (temporarily) blocked during the green phase. The key observation to constructing the mathematical model is that we can model multiple parallel (say $m$) lanes as \emph{one} single queue where batches of (up to) $m$ delayed vehicles can depart in one time slot, for more details see Section~\ref{sec:materials}. The resulting queueing model is one-dimensional just like the standard FCTL queue, which allows us to obtain the probability generating function (PGF) of the steady-state queue-length distribution of the bFCTL queue with multiple lanes and to provide an exact characterization of the capacity. %A slightly different version of the bFCTL queue with a single lane has been studied by means of simulation in a recent paper by Huang et al.~\cite{shaoluen2020random}, which has been the inspiration for the study in the present paper.

	In summary, our main contributions are as follows:
	
	\begin{itemize}
		\item[(i)] We extend the general applicability of the Fixed-Cycle Traffic-Light (FCTL) queue. We allow for traffic streams with multiple lanes and for vehicles to be blocked during the green phase. We refer to this model variation as the blocked Fixed-Cycle Traffic-Light (bFCTL) queue with multiple lanes.
		\item[(ii)] We provide an exact capacity analysis for the bFCTL queue relieving the need for simulation studies.
		\item[(iii)] We provide a way to compute the PGF of the steady-state queue-length distribution of the bFCTL queue and show that it can be used to obtain several performance measures of interest.
		\item[(iv)] We provide a queueing-theoretic framework for the study of shared lanes with potential blockages by pedestrians. This e.g. allows for the study of several performance measures and allows us to model the impact of randomness on the performance measures.
	\end{itemize}
	
	\subsection*{Paper outline}
	The remainder of this paper is organized as follows. In Section~\ref{sec:materials}, we give a detailed model description. This is followed by a capacity analysis, a derivation of the PGF of the steady-state queue-length distribution, and a derivation of some of the main performance measures in Section~\ref{sec:queue_length_derivation}. In Section~\ref{sec:results}, we provide an overview of relevant performance measures for some numerical examples and point out various interesting results. We wrap up with a conclusion and some suggestions for future research in Section~\ref{sec:discussion}.
	
	\section{Detailed model description}\label{sec:materials}
	
	In this section we provide a detailed model description of the bFCTL queue with multiple lanes. %We also highlight the differences with the model described in~\cite{shaoluen2020random} and the traditional FCTL queue.
	
	\begin{figure}[h!]
		\centering
		\begin{tabular}{cc}
			\begin{tikzpicture}[scale=0.1]
				\draw[black,fill=lightgray](30,40) rectangle (40,80);
				\draw[black,fill=lightgray](45,40) rectangle (60,80);
				\draw[black,fill=lightgray](20,40) rectangle (70,50);
				\draw[thick,white](40,40) to (40,80);
				\draw[thick,white](45,40) to (45,80);
				\draw[thick,white,dash pattern=on 7 off 4](35,40) to (35,80);
				\draw[thick,white](50,40) to (50,80);
				\draw[lightgray,fill=lightgray](30,40) rectangle (60,50);
				\draw[black] (40,80) to (50,80);
				\draw[black] (60,50) to (70,50);
				\draw[black] (30,40) to (60,40);
				\filldraw[black] (41,54) circle (6pt);
				\filldraw[black] (42,54) circle (6pt);
				\filldraw[black] (43,54) circle (6pt);
				\filldraw[black] (44,54) circle (6pt);
				
				\draw[green,fill=green](30,50) rectangle (40,50.5);
				\draw[green,fill=green](45,50) rectangle (50,50.5);
				
				\car[270]{(37, 55)}{black}
				\car[270]{(32, 55)}{black}
				\car[270]{(47, 55)}{black}
				\car[270]{(37, 60)}{blue}
				\car[270]{(32, 60)}{blue}
				\car[270]{(47, 60)}{blue}
				\car[270]{(32, 65)}{red}
				\draw [decorate,decoration={brace,amplitude=5pt,mirror,raise=4ex},rotate=180]
				(-29,-55) -- (-51,-55) node[midway,yshift=-3em]{$m$};
				\draw[black] (31,51.5) to (49,51.5);
				\draw[black] (31,51.5) arc(270:90:2.25);
				\draw[lightgray] (31,51.5) -- (31,56);
				\draw[black] (31,56) to (49,56);
				\draw[black] (49,51.5) arc(-90:90:2.25);
				
				\draw[blue] (31,56.5) to (49,56.5);
				\draw[blue] (31,56.5) arc(270:90:2.25);
				\draw[lightgray] (31,56.5) -- (31,61);
				\draw[blue] (31,61) to (49,61);
				\draw[blue] (49,56.5) arc(-90:90:2.25);
				
				\draw[red] (31,61.5) to (49,61.5);
				\draw[red] (31,61.5) arc(270:90:2.25);
				\draw[lightgray] (31,61.5) -- (31,66);
				\draw[red] (31,66) to (49,66);
				\draw[red] (49,61.5) arc(-90:90:2.25);
			\end{tikzpicture} & \hspace{1cm}
			\begin{tikzpicture}[scale=0.10]
				\draw[black,fill=lightgray](0,7) rectangle (10,35);
				\draw[black,fill=lightgray] (5,0) circle (140pt);
				\draw[black,fill=black](1,8) rectangle (9,9);
				\draw[black,fill=black](1,9.5) rectangle (9,10.5);
				\filldraw[black] (5,12.5) circle (3pt);
				\filldraw[black] (5,12) circle (3pt);
				\filldraw[black] (5,11.5) circle (3pt);
				\filldraw[black] (5,11) circle (3pt);
				\draw[black,fill=black](1,13) rectangle (9,14);
				\draw [decorate,decoration={brace,amplitude=5pt,mirror,raise=4ex},rotate=180]
				(-5,-7) -- (-5,-14) node[midway,xshift=3em]{$m$};	
				\draw[black,thick] (0.5,7.5) rectangle (9.5,14.5);
				
				\draw[blue,fill=blue](1,15.5) rectangle (9,16.5);
				\draw[blue,fill=blue](1,17) rectangle (9,18);
				\filldraw[blue] (5,18.5) circle (3pt);
				\filldraw[blue] (5,19) circle (3pt);
				\filldraw[blue] (5,19.5) circle (3pt);
				\filldraw[blue] (5,20) circle (3pt);
				\draw[blue,fill=blue](1,20.5) rectangle (9,21.5);
				\draw [decorate,decoration={brace,amplitude=5pt,mirror,raise=4ex},rotate=180,blue]
				(-5,-15) -- (-5,-22) node[midway,xshift=3em]{$m$};
				\draw[blue,thick] (0.5,15) rectangle (9.5,22);
				
				\draw[red,fill=red](1,23) rectangle (9,24);
				\draw[red,thick] (0.5,22.5) rectangle (9.5,24.5);
			\end{tikzpicture}
			\\
			\scriptsize (a)  %Platoon forming approaching the intersection
			&
			\hspace{0.2cm} \scriptsize (b)  %Crossing the intersection in platoons
		\end{tabular}
		\caption{Visualization of (a) the bFCTL model in terms of an intersection with a traffic stream spread over $m$ lanes and (b) the corresponding queueing model, where the server takes batches of $m$ vehicles into service simultaneously unless there are less than $m$ vehicles present; in that case all vehicles are taken into service.}
		\label{fig:vis_queue}	
	\end{figure}
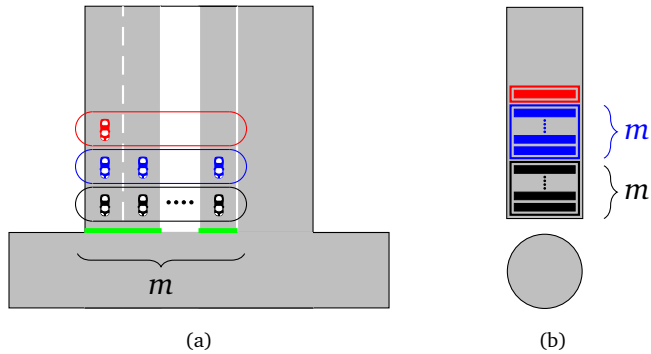
	
	We assume that there are multiple lanes for a traffic stream, that is a group of vehicles coming from the same road and heading into one (or several) direction(s), governed by a \emph{single} traffic light. A visualization can be found in Figure~\ref{fig:vis_queue}(a). As can be seen in Figure~\ref{fig:vis_queue}(a), we assume that there are $m$ lanes and that vehicles spread themselves among the available lanes in such a way that $m$ vehicles can depart if there are at least $m$ vehicles. In practice, this assumption makes sense as drivers gladly minimize their delay by choosing free lanes. The traffic-light model is then turned into a queueing model with a \emph{single} queue with batch services of vehicles, see Figure~\ref{fig:vis_queue}(b). The batches generally consist of $m$ delayed vehicles (we consider delayed vehicles as is done in the study of the FCTL queue, see e.g.~\cite{boon2019pollaczek}), except if less than $m$ delayed vehicles are present at the moment that a batch is taken into service: then all vehicles are taken into service. We further assume that the time axis is divided into time intervals of constant length, where each interval corresponds to the time it takes for a batch of delayed vehicles to depart from the queue. We will refer to these intervals as slots.

	We now turn to discuss two concrete, motivational examples that fit the framework of the bFCTL queue with multiple lanes. After that, we describe the assumptions of the bFCTL queue more formally.
	
	\begin{example}[Shared right-turn lane]
		In this example we consider the scenario as in Figure~\ref{fig:vis}(b). We have batches of vehicles of size $1$, i.e. batches are individual vehicles.
		
		We distinguish between vehicles that are going straight ahead and vehicles that turn right. We do so because only right-turning vehicles can be blocked by crossing pedestrians. The probability that an arbitrary vehicle at the head of the queue is a turning vehicle is $p$. Such a turning vehicle is blocked by a pedestrian in slot $i$ with probability $q_i$, i.e.\,a pedestrian is present on the crossing with probability $q_i$. If a turning vehicle is blocked, all vehicles behind it are also blocked. Then, we proceed to the next slot, $i+1$, and check whether there are any pedestrians crossing (with probability $q_{i+1}$): if there are pedestrians crossing, all vehicles in the queue keep being blocked and otherwise, the turning vehicle at the head of the queue may depart and the blockage of all other vehicles is removed.
		%Then, in the next slot the vehicle at the head of the queue is still a turning vehicle (so when looking at slot $i+1$ we need to know whether or not vehicles during slot $i$ were blocked).
		
		%We further assume that pedestrians only cross during the first part of the green period and not during the second part. This is why we make a division of the green period.
		
		Moreover, if the queue becomes empty during the green period, it will in general not start building again (cf. the FCTL assumption for the regular FCTL queue, see e.g.~\cite{van2006delay}), \emph{except} if there arrives a turning vehicle and there is a crossing pedestrian. The turning vehicle is then blocked and any vehicles arriving in the same slot behind this vehicle are also blocked.
	\end{example}
	
	\begin{example}[Two turning lanes]\label{ex:two}
		In this example we consider the scenario as in Figure~\ref{fig:vis}(c). We have batches of vehicles of size $2$.
		
		In this example, there is no need to make a distinction between vehicles: each vehicle is a turning vehicle with probability $1$, i.e. $p=1$. During each slot $i$, there are pedestrians on the crossing with probability $q_i$ and if there is a pedestrian, all vehicles in the batch are blocked, as are all other vehicles in the queue: there are no vehicles that can complete the right turn. All vehicles in the queue keep being blocked until there are no pedestrians crossing anymore.
		
		Also in this example, the queue of vehicles might dissolve entirely during the green period. If that happens, it only starts building again if there are vehicles arriving \emph{and} if there are pedestrians crossing. In such cases, all arriving vehicles get blocked and remain blocked until there are no pedestrians anymore.
	\end{example}
	
	We are now set to formalize the assumptions for the bFCTL queue with multiple lanes. We number them for clarity and provide additional remarks if necessary. We start with a standard assumption for FCTL queues and a standard assumption on the independence of arriving vehicles, see, e.g.~\cite{van2006delay}.
	
	\begin{assumption} \label{ass:disctime} [Discrete-time assumption]
		We divide time into discrete slots. %!
		The red and green times, $r$ and $g$ respectively, are fixed multiples of those discrete slots and the total cycle length, $c=g+r$, thus consists of an integer number of slots. Each slot corresponds to the duration of the departure of a batch of maximally $m$ delayed vehicles, where $m$ is the maximum number of vehicles that can cross the intersection simultaneously. Any arriving vehicle that finds at least $m$ other vehicles waiting in front of the traffic light is delayed and joins the queue.
	\end{assumption}
	\begin{assumption}[Independence of arrivals]
		All arrivals are assumed to be independent. In particular, the arrivals during slot $i$ do not affect the arrivals in slot $j$ when $i\neq j$.
	\end{assumption}
	
	The next three assumptions, Assumptions~\ref{as:division}, \ref{ass:removal}, and \ref{ass:adapted}, relate to the blockages of vehicles and that allow us to explicitly model such blockages.
	
	\begin{assumption}[Green period division]\label{as:division}
		For the green period we distinguish between two parts, $g_1$ and $g_2$, with $g=g_1+g_2$. During the first part of the green period, blockages might occur (see also Assumption~\ref{ass:removal} below). During the second part of the green period there are no blockages at all. We further assume that $g_2>0$ for technical reasons.
	\end{assumption}

	We make a division of the green period into two parts as is done in e.g.~\cite{shaoluen2020random}. Moreover, such a division is often present in reality and it slightly eases the computations later on. This e.g. means that during the second part of the green period there is a ``no walk'' sign flashing, during which pedestrians are not allowed to cross the intersection. We note that if $g_1=0$ (and $m=1$), we obtain the standard FCTL queue.
	
	Further, we assume that the second part of the green period is strictly positive, mainly for technical reasons. This basically implies that at least one batch of vehicles can depart from the queue during each cycle and that there is \emph{no} batch of vehicles in the queue at the end of the cycle that has caused a blockage before. If $g_2$ would be zero and if a batch of vehicles is blocked at the end of slot $g_1$, this would allow for a blockage to carry over to the next cycle, leading to a slightly more complex model.  Moreover, the red and green times could be taken random in the regular FCTL queue when the times are independent of one another, see e.g.~\cite{boon2021optimal}. At the expense of additional complexity, our framework for the bFCTL queue could be adjusted to account for such sources of randomness. This would allow one to model (to some extent) randomness in, for example, crossing times of pedestrians.
	
	Next, we make an assumption about the blocking of batches of vehicles during the first part of the green period. We take into account that (i) not all batches of vehicles at the head of the queue are potentially blocked (e.g. because only turning batches of vehicles can be blocked); that (ii) if a batch of vehicles is blocked, all vehicles behind it are blocked as well; that (iii) once a blockage occurs, it carries over to the next slot; and that (iv) blockages occur only in the combined event of having a right-turning batch of vehicles at the head of the queue \emph{and} pedestrians crossing the road.
	
	\begin{assumption}[Potential blocking of batches]\label{ass:removal}
		A batch of vehicles, arriving at the head of the queue in time slot $i$, turns right with probability $p_i$. Independently, in time slot $j$, pedestrians cross the road with probability $q_j$, blocking right-turning traffic. As a consequence, whenever a new batch arrives at the head of the queue, this batch will be served in that particular time slot if (i) the batch goes straight ahead, \emph{or} (ii) the batch turns right but there are no crossing pedestrians. Once a batch (of right-turning vehicles) is blocked, it will remain blocked until the next time slot when no pedestrians cross the road. Note that this will be time slot $g_1+1$ at the latest. If the batch at the head of the queue is blocked, it will also block all the other batches in the queue, including those that would go straight. Both $p_i$ and $q_i$ are allowed to depend on the slot $i$.
		
		%A batch that is taken into service during slot $i$ is a right-turning vehicle with probability $p_i$. In slot $i$, If a batch is blocked, all vehicles in that batch get blocked and all vehicles behind it are blocked as well. The vehicles keep being blocked until either slot $g_1+1$ starts (so the first part of the green period is terminated) or until the blockage is no longer in effect.  If $q_i<1$, a blockage might thus not be in effect and even if a batch is susceptible for being blocked, it is \emph{not} blocked in that case. Note that we do not allow for a part of the lanes to be blocked: either all lanes are blocked or none are blocked.
	\end{assumption}
	
	\begin{remark}%[Choices for the $p_i$]
		\label{rem:pi}
				We make a couple of remarks on the values of the $p_i$. 
First, we note that $p_i$ is not representing the probability that the batch at the head of the queue is a turning batch, but rather the probability that a \emph{newly arriving} batch that gets to the head of the queue in slot $i$, is a turning batch. In practice, this will usually \emph{not} depend on the slot in which the batch gets to the head of the queue. This would imply that $p_i=p$ (see, e.g., Example~\ref{ex:two}) and that we could drop the subscript $i$. However, we are able to let $p_i$ depend on the slot in the derivation of the formulas and opt to provide the general case where $p_i$ is allowed to depend on $i$.

		Moreover, in the case that $m>1$, we will often assume that either $p_i=0$, as is the case in Figure~\ref{fig:vis}(a), or $p_i=1$, as is the case in Figure~\ref{fig:vis}(c). This is mainly due to the fact that \emph{all} vehicles in a batch have to be treated similarly: the framework of the bFCTL queue does not allow for batches consisting of one right-turning vehicle that is blocked and one straight-going vehicle that is allowed to depart because it is not blocked. I.e. a case with mixed traffic and \emph{multiple} lanes, such as the shared right-turn lane example in Figure~\ref{fig:vis}(b) but with $m>1$, is not modeled by the bFCTL queue. We do not consider this to be a severe restriction as it will often be the case in practice that $p_i=0$ or $p_i=1$ if $m>1$. We stress that the case with $m=1$ as depicted by the blue rectangle in Figure~\ref{fig:vis}(b) can be studied by the bFCTL queue.
	\end{remark}
	
	\begin{remark}%[Blockages extend over several slots]
		\label{rem:blockages}
		We would like to stress that the blockage of a batch of vehicles carries over to the next slot. E.g. if a vehicle is a right-turning vehicle in Figure~\ref{fig:vis}(b) and is blocked, it is still at the head of the queue in the next slot. So, as soon as a blockage actually takes place, we are essentially in a different state of the system than in the case where there is no blockage: if there is a blockage in time slot $i$ then we are sure that there is a right-turning batch at the head of the queue in time slot $i+1$. This is why we have two mechanisms for the blocking: on the one hand we have the $p_i$ to check whether \emph{new} batches that get to the head of the queue are right turning and on the other hand we have the pedestrians crossing in slot $i$ accounted for by the $q_i$.
	\end{remark}
	
	%Assumption~\ref{ass:removal} generalizes~\cite{shaoluen2020random} where $q_i$ is chosen to be $1$ for $i=1,\dots,g_1$. The possibility to choose $q_i<1$ allows us to model an intersection with only a few pedestrians or where the expected number of crossing pedestrians decreases during the first part of the green period. Note that the fact that the $p_i$ and $q_i$ depend on slot $i$ implies that there is no strict need to distinguish between the first and second part of the green period as in Assumption~\ref{as:division}. Indeed, e.g. by choosing the $q_i$ properly (i.e. $q_i=0$ for any green slot in the second part of the green phase), we create the exact same effect. For convenience, we keep referring to the first and second part of the green phase, as is done in~\cite{shaoluen2020random}, and because in practice such a division is often present.
	
	We need one final assumption which is a slightly adapted version of the standard FCTL assumption. We require a slight change because of the potential blocking of vehicles during the first part of the green phase and because of the possibility that there is more than one delayed vehicle departing in a single slot during the green period because of the batch-service structure.
	\begin{assumption}[bFCTL assumption]\label{ass:adapted}
		We assume that any vehicle arriving during a slot where $m-1$ or less vehicles are in the queue, may depart from the queue immediately together with the $m-1$ or less delayed vehicles. There are two exceptions: (i) if this batch of $m-1$ or less vehicles is blocked or (ii) if the queue was empty and there is an arriving vehicle that gets blocked, in which case that vehicle gets blocked together with any arriving vehicles after that vehicle. In the former case, all arriving vehicles together with the delayed vehicles remain at the queue. In the latter case, the first blocked vehicle is delayed and any arriving vehicles behind it (if any) are also delayed and blocked where we restrict ourselves to the situation where the queue is empty. If the queue was not empty, then we assume that either all arriving vehicles in that slot are blocked and delayed (because the batch at the head of the queue is blocked) or that all arriving vehicles are allowed to depart along with the batch of delayed vehicles (because the batch at the head of the queue is not blocked). Summarizing, if the queue length at the start of the slot is at least $1$ but at most $m-1$, we either have no departures (in case of a blockage) or \emph{all} vehicles are allowed to cross the intersection (including arriving vehicles). If the queue length is $0$, we only have a non-zero queue at the end of the slot if a vehicle gets blocked: then the blocked vehicle and any vehicles arriving behind it are queued.
	\end{assumption}
	
	\begin{remark}%[Concrete examples for the bFCTL queue]
		\label{rem:bFCTL}
		The bFCTL assumption allows one to model a situation where arriving vehicles get blocked if the queue was already empty before the start of the slot. Although, in principle, one can use any distribution for the number of arriving vehicles that are blocked, there are only few logical choices in practice. For example, in the case of Figure~\ref{fig:vis}(b), the number of (potentially) blocked vehicles that arrive at the queue during slot $i$ would correspond to the number of vehicles counting from the first right-turning vehicle among all vehicles arriving in slot $i$: these vehicles will be blocked if there is a crossing pedestrian in slot $i$. In Figure~\ref{fig:vis}(c), any arriving vehicle is a turning vehicle. So, if there is a crossing  pedestrian, all arriving vehicles in slot $i$ are blocked.
	\end{remark}
	
	The combination of all the above assumptions enables us to view the process as a discrete-time Markov chain, which in turn allows us to obtain the capacity and the PGF of the steady-state queue-length distribution of the bFCTL queue with multiple lanes. We do so in the next section.
	
	\section{Capacity analysis, PGFs, and performance measures for the bFCTL queue}\label{sec:queue_length_derivation}
	
	In this section we provide an exact analysis for the bFCTL queue. We start with an exact characterization of the capacity in Subsection~\ref{subsec:capacity}. In Subsection~\ref{subsec:pgf_derivation}, we obtain the steady-state queue-length distribution in terms of PGFs where we thus focus on the \emph{transforms} of the queue-length distribution, because we cannot directly obtain closed-form expressions for the probabilities. We can use the methods devised in e.g.~\cite{abate1992numerical} and~\cite{choudhury1996numerical} to obtain numerical values from the PGFs for the queue-length probabilities and moments respectively. Without giving details, we stress that our recursive approach in Subsection~\ref{subsec:pgf_derivation} also allows us to provide a transient analysis in which case we can also take time-varying parameters into account. In Subsection~\ref{sec:performance_measures}, we study several important performance measures of the bFCTL queue.
	
	\subsection{Capacity analysis for the bFCTL queue}\label{subsec:capacity}
	
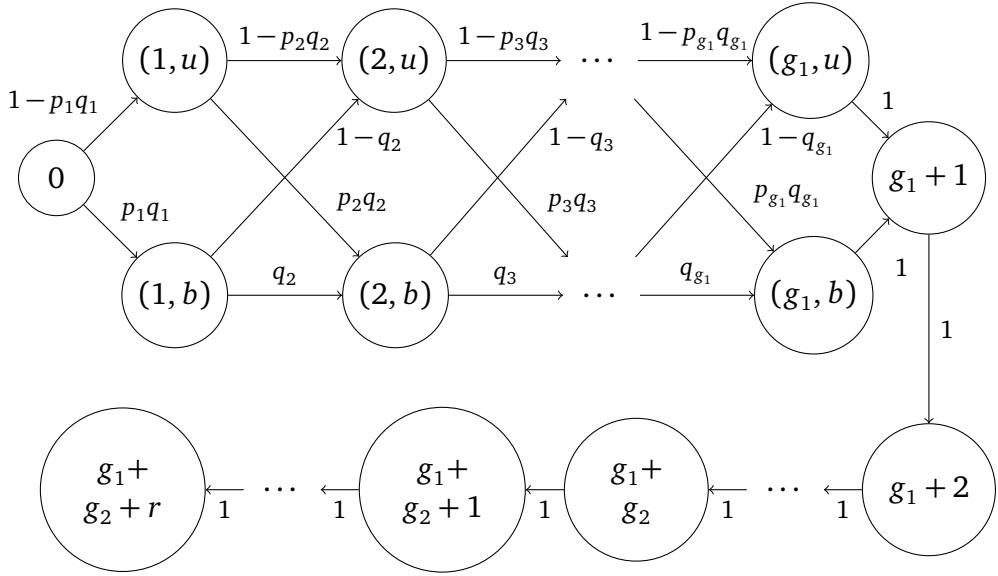
\begin{figure}[!ht]
		\centering
		\begin{tikzpicture}[->,auto,node distance=2.2cm]
			\tikzstyle{every state}=[fill=white,draw=black,text=black,minimum size=1.1cm]
			\node[draw,circle,minimum height = 1cm] (A) {$0$};
			\node[draw,circle,minimum height = 1cm] (B) [above right of=A] {$(1,u)$};
			\node[draw,circle,minimum height = 1cm] (C) [below right of=A] {$(1,b)$};
			\node[draw,circle,minimum height = 1cm] (D) [right  = 1.5 cm of B] {$(2,u)$};
			\node[draw,circle,minimum height = 1cm] (E) [right  = 1.5 cm of C] {$(2,b)$};
			\node[minimum size=1cm] (H) [right=1.5 cm of D] {$\dots$};
			\node[minimum size=1cm] (I) [right=1.5 cm of E] {$\dots$};
			\node[draw,circle,minimum height = 1cm] (J) [right=1.5 cm of H] {$(g_1,u)$};
			\node[draw,circle,minimum height = 1cm] (K) [right=1.5 cm of I] {$(g_1,b)$};
			\node[draw,circle,minimum height = 1cm] (L) [below right of=J] {$g_1+1$};
			\node[draw,circle,minimum height = 1.75cm] (LL) [below=2.5 cm of L] {$g_1+2$};
			\node[minimum size=1cm] (M) [left=0.5 cm of LL] {$\dots$};	
			\node[draw,circle,minimum height = 1.75cm] (N) [left=0.5 cm of M] {\begin{tabular}{c}$g_1+$\\$g_2$\end{tabular}};
			\node[draw,circle,minimum height = 2cm] (O) [left=0.5 cm of N] {\begin{tabular}{c}$g_1+$\\$g_2+1$\end{tabular}};
			\node[minimum size=1cm] (P) [left=0.5 cm of O] {$\dots$};	
			\node[draw,circle,minimum height = 1.75cm] (Q) [left=0.5 cm of P] {\begin{tabular}{c}$g_1+$\\$g_2+r$\end{tabular}};
			
			\path (A) edge node {\footnotesize{$1-p_1q_1$}} (B);
			\path (A) edge node {\footnotesize{$p_1q_1$}} (C);
			\path (B) edge node {\footnotesize{$1-p_2q_2$}} (D);
			\path (B) edge node[pos=0.78] {\footnotesize{$p_2q_2$}} (E);
			\path (C) edge node[pos=0.6,xshift=1.5cm] {\footnotesize{$1-q_2$}} (D);
			\path (C) edge node {\footnotesize{$q_2$}} (E);
			\path (D) edge node {\footnotesize{$1-p_3q_3$}} (H);
			\path (D) edge node[pos=0.78] {\footnotesize{$p_3q_3$}} (I);
			\path (E) edge node[pos=0.6,xshift=1.5cm] {\footnotesize{$1-q_3$}} (H);
			\path (E) edge node {\footnotesize{$q_3$}} (I);
			\path (H) edge node {\footnotesize{$1-p_{g_1}q_{g_1}$}} (J);
			\path (H) edge node[pos=0.78] {\footnotesize{$p_{g_1}q_{g_1}$}} (K);
			\path (I) edge node[pos=0.6,xshift=1.7cm] {\footnotesize{$1-q_{g_1}$}} (J);
			\path (I) edge node {\footnotesize{$q_{g_1}$}} (K);
			\path (K) edge node[pos=0.1,xshift=0.8cm,yshift=-0.5cm] {\footnotesize{$1$}} (L);
			\path (J) edge node {\footnotesize{$1$}} (L);
			\path (L) edge node {\footnotesize{$1$}} (LL);
			\path (LL) edge node {\footnotesize{$1$}} (M);
			\path (M) edge node {\footnotesize{$1$}} (N);
			\path (N) edge node {\footnotesize{$1$}} (O);
			\path (O) edge node {\footnotesize{$1$}} (P);
			\path (P) edge node {\footnotesize{$1$}} (Q);
		\end{tikzpicture}
		\caption{Markov chain used to study the capacity of the bFCTL queue.}\label{fig:MC}
	\end{figure}

	In this subsection we develop a computational algorithm to determine the capacity for the bFCTL queue. The capacity is defined as the maximum  number of vehicles that can cross the intersection in the given lane group, per time unit. In the standard FCTL queue, the capacity can simply be determined by multiplying the saturation flow with the ratio of the green time and the cycle length. In the bFCTL model, however, there are subtle dependencies which carry over from one cycle to the next cycle. We will capture these dependencies by means of a Markov reward model. The Markov chain with the associated transition probabilities that we use is depicted in Figure~\ref{fig:MC}. We are interested in the number of departures of delayed vehicles in each time slot. For this reason, the Markov chain that we consider here only has states $(i,s)$ for $i=1,\dots,g_1$ representing the slots during the first part of the green period and $s=u,b$, representing the case where vehicles are not blocked ($s=u$) and the case where vehicles are blocked ($s=b$). We also have states $i$ for $i=g_1+1,\dots,g_1+g_2+r$ representing the slots during the second part of the green period and the red period. Finally, we create an artificial state $0$ to gather the rewards from states $(1,b)$ and $(1,u)$. The long-term mean number of departures of delayed vehicles can now be determined by means of a Markov reward analysis.
	
	We use Markov reward theory to obtain the capacity of the bFCTL queue. In order to use Markov reward theory, we work backwards from state $g_1+g_2+r$ to obtain the reward in state $0$. Indeed, we get the mean number of vehicles that is able to depart from the queue in an arbitrary cycle when we compute the reward in state $0$. The rewards that we assign to each transition are as follows: if we make a transition to a state $(i,u)$ for $i=1,\dots,g_1$, we receive a reward $m$ reflecting the maximum of $m$ delayed vehicles departing from the queue. We also get a reward $m$ if we make a transition from state $g_1+i$ to state $g_1+i+1$ for $i=1,\dots,g_2-1$. For all other transitions, we receive no reward as there are no vehicles departing. We denote the received reward up to state $(i,s)$ with $r_{i,s}$ with $i=1,\dots,g_1$ and $s=u,b$ and the received reward up to state $i$ with $r_i$ for $i=0$ and $i=g_1+1,\dots,g_1+g_2+r$. Then we get the following relations between the rewards in the various states. We start with defining the total reward in state $g_1+g_2+r$ to be $0$ (there are no vehicle departures while being in state $g_1+g_2+r$), i.e.
	\begin{equation} \label{eq:rc}
		r_{g_1+g_2+r} = 0.
	\end{equation}
	For states $i=g_1+g_2,\dots,g_1+g_2+r-1$, we obtain
	\begin{equation}
		r_{i} = r_{i+1},
	\end{equation}
	as there are no departures during the red period. However, for states $i=g_1+1,\dots,g_1+g_2-1$, we have
	\begin{equation}\label{eqnm1}
		r_{i} = m + r_{i+1}
	\end{equation}
	as there are (potentially) $m$ delayed vehicles departing. For state $(g_1,b)$ we have that
	\begin{equation}
		r_{g_1,b} = r_{g_1+1},
	\end{equation}
	as there are no departures when the vehicles are blocked. For state $(g_1,u)$ we obtain
	\begin{equation}
		r_{g_1,u} = m + r_{g_1+1}\label{eqnm2}
	\end{equation}
	as there are, at most, $m$ delayed vehicles departing from the queue when we transition from state $(g_1,u)$ to $g_1+1$. Similarly, for states $(i,b)$ with $i=1,\dots,g_1-1$, we get
	\begin{equation}
		r_{i,b} = q_{i+1}r_{i+1,b}+(1-q_{i+1})r_{i+1,u}
	\end{equation}
	and for states $(i,u)$ with $i=1,\dots,g_1-1$, we get
	\begin{equation}\label{eqnm3}
		r_{i,u} = m+p_{i+1}q_{i+1}r_{i+1,b}+(1-p_{i+1}q_{i+1})r_{i+1,u}.
	\end{equation}
	Finally, for state $0$, we get
	\begin{equation}\label{eq:r0}
		r_{0} = p_1q_1r_{1,b}+(1-p_1q_1)r_{1,u}.
	\end{equation}
	
	Then we have that $r_0$ is the average reward received when traversing the Markov chain as depicted in Figure~\ref{fig:MC}. This average reward translates to the mean number of delayed vehicles that are able to depart from the queue during a cycle, which is exactly the capacity of this lane group. We can thus compute the capacity of the bFCTL queue for each set of input parameters. Along with the mean number arrivals per cycle, we can also check whether the bFCTL queue renders a stable queueing model. If we denote the mean number of arrivals in slot $i$ by $\mathbb{E}[Y_i]$, the mean number of arrivals per cycle is $\sum_{i=1}^c\mathbb{E}[Y_i]$ and the bFCTL queue is stable if $r_0>\sum_{i=1}^c\mathbb{E}[Y_i]$. The procedure to check for stability is summarized in Algorithm~\ref{alg:stability}.
	
	\begin{algorithm}[H]
		\caption{Algorithm to check for stability of the bFCTL queue.} \label{alg:stability}
		\begin{algorithmic}[1]
			\State Input: $\mathbb{E}[Y_i]$ for $i=1,\dots,c$, $g_1$, $g_2$, $c$, $p_{i}$ for $i=1,\dots,g_1$, and $q_{i}$ for $i=1,\dots,g_1$.
			\State Use Equations~\eqref{eq:rc} up to \eqref{eq:r0} to determine $r_0$.
			\If {$\sum_{i=1}^c\mathbb{E}[Y_i]< r_0$}
			\State The bFCTL queue is stable.
			\Else
			\State The bFCTL queue is not stable.
			\EndIf
		\end{algorithmic}
	\end{algorithm}
	
\begin{remark}\label{rem:capacity}
One of our model restrictions (Assumption~\ref{ass:disctime}) is that vehicles depart at the end of each time slot, meaning that we do not correct for the fact that turning vehicles might need more time to accelerate. A simple method to account for this effect, which reduces the capacity in practice, is to modify the reward structure of the Markov chain. One can modify the value of $m$ in Equations~\eqref{eqnm1}, \eqref{eqnm2}, and \eqref{eqnm3} to account for the lower departure rate of turning vehicles. For example, one can use
\begin{equation}\label{eqn:capacityCorrection}
m^* = p_i m_\textit{turn} + (1-p_i)m_\textit{through},
\end{equation}
where $m_\textit{through}$ and $m_\textit{turn}$ represent the average number of through-vehicles and turning vehicles, respectively, crossing the intersection per time unit. For this capacity calculation, these numbers do not need to be integers. See Section~\ref{subsec:capacity2} for a numerical example and a comparison to the HCM capacity formula.
\end{remark}

	\subsection{Derivation of the PGFs for the bFCTL queue}\label{subsec:pgf_derivation}
	
	First, we need to introduce some further concepts and notation before we continue our quest to obtain the relevant PGFs of the queue-length distribution. We introduce two states, one corresponding to a situation where the queue is blocked and one where this is not the case, cf. Assumption~\ref{ass:removal} and Remark~\ref{rem:blockages} and as is done in Subsection~\ref{subsec:capacity}. We denote the random variable of being in either of the two states with $S$ and $S$ takes the values $b$ (blocked) and $u$ (unblocked). By definition, blocked states only occur during the first part of the green period and if there are vehicles in the queue. We define $S$ to be equal to $u$ if the queue is empty. We denote the joint steady-state queue length (measured in number of vehicles) and the state $S$ at the end of slot $i=1,\dots,g_1$ with the tuple $(X_i,S)$ and we denote its PGF with $X_{i,j}(z)$ where $i=1,\dots,g_1$ and $j=u,b$. We note that $X_{i,b}(z)$ and $X_{i,u}(z)$ are partial generating functions: we e.g. have $X_{i,b}(z) = \mathbb{E}[z^{X_{i}}\mathds{1}\{S=b\}]$, where $\mathds{1}\{S=b\}=1$ if $S=b$ and $0$ otherwise. For the slots $i=1,\dots,c$ we denote the steady-state queue length with $X_i$ and its PGF with $X_i(z)$, so for $i=1,\dots,g_1$ we have that $X_i(z)=X_{i,u}(z)+X_{i,b}(z)$.
	
	We note that, as we are looking at the steady-state distribution of the number of vehicles in the queue, we need to require stability of the queueing model. We can check whether or not the stability condition is satisfied by means of Algorithm~\ref{alg:stability} devised in Subsection~\ref{subsec:capacity}.
	
	We further denote with $Y_i$ the number of arrivals during slot $i$ and with $Y_{i,b}$ we denote the total number of arrivals of potentially blocked vehicles during slot $i$, see also Assumption~\ref{ass:adapted}. %We only have to consider $Y_{i,b}$ if the queue was empty at the start of the slot: if there were $m$ or more vehicles at the start of the slot, all arriving vehicles are automatically delayed, whereas if there were $m-1$ vehicles or less, either ($1$) all arriving vehicles in the slot are delayed if the batch of at most $m-1$ vehicles is blocked or ($2$) the batch was not blocked, in which case
	We denote their PGFs respectively with $Y_i(z)$ and $Y_{i,b}(z)$. Later in this subsection, we provide $Y_{i,b}(z)$ for several concrete examples.

	In the next part of this subsection we provide the recursion between the $X_{i,j}(z)$, $i=1,\dots,g_1$ and $j=u,b$, and the $X_i(z)$, $i=g_1+1,\dots,c$. Afterwards, we wrap up with some technicalities that need to be overcome to obtain a full characterization of all the PGFs.
	
	\subsubsection{Recursion for the $X_{i,j}(z)$}
	
	We start with the relation between $X_{1,b}(z)$ and $X_c(z)$. We distinguish several cases while making a transition from slot $c$ to a blocked state in slot $1$. We get
	\begin{equation}\label{eq:X1b}
		\begin{aligned}
			X_{1,b}(z) = & p_1q_1 \mathbb{E}[z^{X_c+Y_1}\mathds{1}\{X_c>0\}] + q_1 \mathbb{E}[z^{Y_{1,b}}\mathds{1}\{X_c=0\}\mathds{1}\{Y_{1,b}>0\}] +\\& 0\cdot\mathbb{E}[ \mathds{1}\{X_c=0\}\mathds{1}\{Y_{1,b}=0\}]\\
			= & p_1q_1 X_{c}(z) Y_1(z) + q_1\mathbb{P}(X_{c}=0)\left(Y_{1,b}(z)-Y_{1,b}(0)-p_1Y_1(z)\right).
		\end{aligned}
	\end{equation}
	We explain this relation as follows: if the queue is nonempty at the end of slot $c$, we need both a right-turning batch of vehicles and a crossing pedestrian in slot $1$ to get a blockage, which happens with probability $p_1q_1$. The queue length at the end of slot $1$ is then $X_c+Y_1$. The second term can be understood as follows: if $X_c=0$, the queue at the end of slot $c$ is empty and then we get to a blocked state if there is a pedestrian crossing (which happens with probability $q_1$) and if $Y_{1,b}>0$, in which case the queue length is $Y_{1,b}$. Note that we further have that the case $X_{1,b}=0$ cannot occur (by definition) as indicated by the term on the second line of Equation~\eqref{eq:X1b}.
	
	Similarly, we derive $X_{1,u}(z)$:
	\begin{align}\label{eq:X1u}
		X_{1,u}(z)  = &\nonumber (1-p_1q_1) \mathbb{E}[z^{X_c+Y_1-m}\mathds{1}\{X_c\geq m\}] + (1-p_1q_1)\mathbb{E}[z^0\mathds{1}\{1\leq X_c\leq m-1\}] + \\&(1-q_1)\mathbb{E}[z^0 \mathds{1}\{X_c=0\}]+q_1\mathbb{E}[z^0 \mathds{1}\{X_c=0\}\mathds{1}\{Y_{1,b}=0\}] \\
		= &\nonumber(1-p_1q_1) X_{c}(z) \frac{Y_1(z)}{z^m} + (1-p_1q_1)\sum_{l=1}^{m-1}\mathbb{P}(X_{c}=l)\left(1-\frac{Y_{1}(z)}{z^{m-l}}\right)+\\&\mathbb{P}(X_c=0)\left(1-q_1+q_1Y_{1,b}(0)-(1-p_1q_1)\frac{Y_1(z)}{z^m}\right).\nonumber
	\end{align}
	This relation can be understood in the following way: first, if there are at least $m$ vehicles at the end of slot $c$ and if there is no blockage (which occurs with probability $1-p_1q_1$, i.e. the complement of a blockage occurring), then the queue length at the end of slot $1$ is $X_c+Y_1-m$. Secondly, if there is at least $1$ but at most $m-1$ vehicles at the end of slot $c$, we have an empty queue at the end of slot $1$ if there is no blockage (which is the case with probability $1-p_1q_1$). Thirdly, if the queue is empty at the end of slot $c$, then the queue remains empty if there are no pedestrians crossing (occurring with probability $1-q_1$) or if there is a pedestrian crossing (occurring with probability $q_1$) while $Y_{1,b}=0$. This fully explains Equation~\eqref{eq:X1u}.
	
	In a similar way, we obtain the following relations for slots $i=2,\dots,g_1$:
	\begin{equation}\label{eq:Xib}
		\begin{aligned}
			X_{i,b}(z)  = & p_iq_i \mathbb{E}[z^{X_{i-1}+Y_i}\mathds{1}\{S=u\}] + q_i \mathbb{E}[z^{X_{i-1}+Y_i}\mathds{1}\{S=b\}] + \\& q_i\mathbb{E}[z^{Y_{i,b}}\mathds{1}\{X_{i-1}=0\}\mathds{1}\{S=u\}\mathds{1}\{Y_{i,b}>0\}]\\
			= &p_iq_i X_{i-1,u}(z) Y_i(z) + q_iX_{i-1,b}(z)Y_i(z) +\\&q_i\mathbb{P}(X_{i-1}=0,S=u)\left(Y_{i,b}(z)-Y_{i,b}(0)-p_iY_{i}(z)\right),
		\end{aligned}
	\end{equation}
	where we have to take both transitions from slot $i-1$ while being blocked (the case $S=b$) and not being blocked (the case $S=u$) into account, and
	\begin{align}
		X_{i,u}(z)  = & (1-p_iq_i) \mathbb{E}[z^{X_{i-1}+Y_i-m}\mathds{1}\{X_{i-1}\geq m\}\mathds{1}\{S=u\}] +\nonumber\\& (1-q_i) \mathbb{E}[z^{X_{i-1}+Y_i-m}\mathds{1}\{X_{i-1}\geq m\}\mathds{1}\{S=b\}]+ \nonumber\\& (1-p_iq_i)\mathbb{E}[z^0\mathds{1}\{1\leq X_{i-1}\leq m-1\}\mathds{1}\{S=u\}] + \nonumber\\& (1-q_i)\mathbb{E}[z^0\mathds{1}\{1\leq X_{i-1}\leq m-1\}\mathds{1}\{S=b\}] +\nonumber\\& (1-q_i)\mathbb{E}[z^0\mathds{1}\{X_{i-1}=0\}\mathds{1}\{S=u\}]+\nonumber\\&q_i\mathbb{E}[z^0\mathds{1}\{X_{i-1}=0\}\mathds{1}\{S=u\}\mathds{1}\{Y_{i-1,b}=0\}] \\
		= &(1-p_iq_i) X_{i-1,u}(z) \frac{Y_i(z)}{z^m} + (1-q_i)X_{i-1,b}(z)\frac{Y_i(z)}{z^m}+ \nonumber\\&(1-p_iq_i)\sum_{l=1}^{m-1}\mathbb{P}(X_{i-1}=l,S=u)\left(1-\frac{Y_{i}(z)}{z^{m-l}}\right)+\nonumber\\& (1-q_i)\sum_{l=1}^{m-1}\mathbb{P}(X_{i-1}=l,S=b)\left(1-\frac{Y_{i}(z)}{z^{m-l}}\right)+\nonumber\\&\mathbb{P}(X_{i-1}=0,S=u)\left(1-q_i+q_iY_{i,b}(0)-(1-p_iq_i)\frac{Y_i(z)}{z^m}\right).\nonumber
	\end{align}
	In order to derive $X_{g_1+1}(z)$, we note that we need to take the cases into account where the queue was blocked or not during slot $g_1$. We then get
	\begin{align}
		\begin{aligned}
			X_{g_1+1}(z) = & \mathbb{E}[z^{X_{g_1}+Y_{g_1+1}-m}\mathds{1}\{X_{g_1}\geq m\}\mathds{1}\{S=u\}]+\\&\mathbb{E}[z^{X_{g_1}+Y_{g_1+1}-m}\mathds{1}\{X_{g_1}\geq m\}\mathds{1}\{S=b\}]+\\&\mathbb{E}[z^0\mathds{1}\{X_{g_1}\leq m-1\}\mathds{1}\{S=u\}]+\mathbb{E}[z^0\mathds{1}\{X_{g_1}\leq m-1\}\mathds{1}\{S=b\}] \\
			= & \left(X_{g_1,u}(z) +X_{g_1,b}(z) \right)\frac{Y_{g_1+1}(z)}{z^m}+\\&\sum_{l=0}^{m-1}\left(\mathbb{P}(X_{g_1}=l,S=u\right)\left(1-\frac{Y_{g_1+1}(z)}{z^{m-l}}\right)+\\&\sum_{l=1}^{m-1}\left(\mathbb{P}(X_{g_1}=l,S=b\right)\left(1-\frac{Y_{g_1+1}(z)}{z^{m-l}}\right).
		\end{aligned}
	\end{align}
	For $i=g_1+2,\dots,g_1+g_2$, we obtain the following
	\begin{equation}
		\begin{aligned}
			X_{i}(z) = & \mathbb{E}[z^{X_{i-1}+Y_{i}-m}\mathds{1}\{X_{i-1}\geq m\}]+\mathbb{E}[z^0\mathds{1}\{X_{i-1}\leq m-1\}] \\
			= & X_{i-1}(z) \frac{Y_{i}(z)}{z^m}+\sum_{l=0}^{m-1}\left(\mathbb{P}(X_{i-1}=l\right)\left(1-\frac{Y_{i}(z)}{z^{m-l}}\right),
		\end{aligned}
	\end{equation}
	while for slots $i=g_1+g_2+1,\dots,c$ we get
	\begin{equation}\label{eq:Xg1+g2+i}
		\begin{aligned}
			X_{i}(z) = \mathbb{E}[z^{X_{i-1}+Y_{i}}]= X_{i-1}(z)Y_i(z).
		\end{aligned}
	\end{equation}

	The combination of all equations above, provides us with a recursion with which we can express $X_{g_1+g_2}(z)$ in terms of $Y_i(z)$, $Y_{i,b}(z)$, $\mathbb{P}(X_{i}=l,S=u)$ and $\mathbb{P}(X_{i}=l,S=b)$ for $i=1,\dots,g_1$ and $l=0,\dots,m-1$, and $\mathbb{P}(X_{i} = l)$ for $i=g_1+1,\dots,g_1+g_2-1$, $i=c$, and $l=0,\dots,m-1$, with the following general form:
	\begin{equation}\label{eq:generalbFCTL}
		X_{g_1+g_2}(z) = \frac{X_n(z)}{X_d(z)},
	\end{equation}
	with known $X_n(z)$ and $X_d(z)$. We refrain from giving $X_n(z)$ and $X_d(z)$ in the general case because of their complexity and only provide them under simplifying assumptions later in this subsection. The $Y_{i}(z)$ are known, but we still need to obtain the $Y_{i,b}(z)$, the $\mathbb{P}(X_{i}=l,S=u)$ and $\mathbb{P}(X_{i}=l,S=b)$ for $i=1,\dots,g_1$ and $l=0,\dots,m-1$, and the $\mathbb{P}(X_{i} = l)$ for $i=g_1+1,\dots,g_1+g_2-1$, $i=c$, and $l=0,\dots,m-1$. We start with the $Y_{i,b}(z)$ and then come back to the unknown probabilities.
	
	The occurrence of the PGF $Y_{i,b}(z)$ directly relates to Assumption~\ref{ass:adapted}. As mentioned before in Remark~\ref{rem:bFCTL}, one could, a priori, use any positively distributed, discrete random variable. However, when we have a specific example in mind, there is usually one logical definition, see also Remark~\ref{rem:Yib} below.
	
	\begin{remark}\label{rem:Yib}%[The bFCTL assumption revisited]
		In general, we define $Y_{i,b}$ to be the random variable of the total number of arrivals of potentially blocked vehicles during slot $i$, cf. Assumption~\ref{ass:adapted}. In case $m=1$, such as in Figure~\ref{fig:vis}(b), the interpretation of the $Y_{i,b}(z)$ is straightforward. We simply count the number of arriving vehicles starting from the first vehicle that is a turning vehicle. We get the following expression for $Y_{i,b}(z)$:
		\begin{align*}
			Y_{i,b}(z) & = \sum_{k=0}^\infty \mathbb{P}(Y_{i,b} = k)z^k\\
			& = \sum_{j=0}^\infty \mathbb{P}(Y_i=j)(1-p_i)^j + \sum_{k=1}^\infty \sum_{j=k}^\infty \mathbb{P}(Y_i = j) (1-p_i)^{j-k} p_i z^k \\
			& = Y_i(1-p_i) + \sum_{j=1}^\infty p_i \mathbb{P}(Y_i=j)(1-p_i)^j\sum_{k=1}^j \left(\frac{z}{1-p_i}\right)^k \\
			& = Y_i(1-p_i) + \sum_{j=1}^\infty p_i \mathbb{P}(Y_i=j)(1-p_i)^j z\frac{1-\left(\frac{z}{1-p_i}\right)^j}{1-p_i-z}\\
			& = Y_i(1-p_i) + \frac{p_i z}{1-p_i-z}\sum_{j=1}^\infty  \mathbb{P}(Y_i=j)\left((1-p_i)^j-z^j\right)\\
			%& = Y_i(1-p_i) + \frac{p_iz}{1-p_i-z} \left(Y_i(1-p_i)-Y_i(0) - Y_i(z) + Y_i(0)\right) \\
			& = Y_i(1-p_i) + \frac{p_iz}{1-p_i-z}\left(Y_i(1-p_i)-Y_i(z)\right),
		\end{align*}
		where in the second step we condition on the total number of arrivals and take into account how we can get to $k$ blocked vehicles; in the third step we interchange the order of the summation; and in the fourth step we compute a geometric series.
		
		If $m>1$, the interpretation as above for the case $m=1$ is not necessarily meaningful. It is more difficult to compute the $Y_{i,b}$ in a logical and consistent way. This has to do with the fact that if $m>1$ we consider batches of vehicles that are either all blocked or not, whereas the $Y_{i,b}$'s are about individual vehicles. As mentioned before in Remark~\ref{rem:pi}, if $m>1$ we often have that either $p_i=0$ or $p_i=1$. If $p_i=0$, the general expression for $Y_{i,b}(z)$ reduces to:
		\[
		Y_{i,b}(z) = Y_{i}(1) + 0\cdot(Y_i(1)-Y_i(z)) = Y_{i}(1)=1,
		\]
		which makes sense as there are no turning vehicles in case $p_i=0$. If $p_i=1$, we have that:
		\[
		Y_{i,b}(z) = Y_{i}(0) - (Y_{i}(0)-Y_{i}(z))=Y_{i}(z),
		\]
		which is also logical: every arriving vehicle is a turning vehicle if $p_i=1$, so we have that $Y_{i,b}(z)=Y_{i}(z)$.
	\end{remark}
	
	Except for the constants $\mathbb{P}(X_{i}=l,S=u)$ and $\mathbb{P}(X_{i}=l,S=b)$ for $i=1,\dots,g_1$ and $l=0,\dots,m-1$, and $\mathbb{P}(X_{i} = l)$ for $i=g_1+1,\dots,g_1+g_2-1$, $i=c$, and $l=0,\dots,m-1$, we are now done. We explain how to find the (so far) unknown constants in the next part of this subsection.
	\subsubsection{Finding the unknowns in $X_{g_1+g_2}(z)$}\label{subsubsec:completion}
	
	As mentioned before, we still need to find several unknowns before the expression for $X_{g_1+g_2}(z)$ is complete. %We employ the same arguments to The denominator in Equation~\eqref{eq:FCTL} has $g$ zeros in the complex plane which satisfy $|z|\leq 1$ whenever the queue is operating below maximum capacity (i.e. $g < c \mathbb{E}[Y]$). This can be shown using Rouch\'{e}'s theorem, see e.g.~\cite{adan2006application} . As $X_{g}(z)$ is a PGF and therefore analytic for all $|z|\leq 1$, we need the numerator in Equation~\eqref{eq:FCTL} to vanish at all zeros of the denominator with a modulus smaller than or equal to $1$. One of the zeros of the denominator is $z=1$, which also leads to a zero of the numerator regardless of the unknown constants. We thus need one more equation and we find this in the requirement that for $z=1$ any PGF should be equal to one. Using L'H\^{o}pital's rule, we can find this normalization equation. Then, we have a set of $g$ linear equations with $g$ unknowns, which we might solve to find the unknown terms $X_{i}(0)=\mathbb{P}(X_{i}=0)$. This system is explicitly captured in Equation~(12) of~\cite{van2006delay} and turns out to be a Vandermonde system.
	The standard framework for the FCTL queue as described in e.g.~\cite{van2006delay} is also applicable to the bFCTL queue with multiple lanes with some minor differences. Although we are dealing with more complex formulas, the key ideas are identical. We have $m(g_1+g_2)+(m-1)g_1$ unknowns in the numerator $X_n(z)$ of $X_{g_1+g_2}(z)$  in Equation~\eqref{eq:generalbFCTL} and we have $m(g_1+g_2)$ roots with $|z|\leq 1$ for the denominator $X_d(z)$ of $X_{g_1+g_2}(z)$, assuming stability of the queueing model. An application of Rouch\'{e}'s theorem, see e.g.~\cite{adan2006application}, shows that $X_d(z)$ indeed has $m(g_1+g_2)$ roots on or within the unit circle assuming stability. One root is $z=1$, which leads to a trivial equation and as a substitute for this root, we put in the additional requirement that $X_{g_1+g_2}(1)=1$. The remaining $(m-1)g_1$ equations are implicitly given in Equations~\eqref{eq:X1b} and \eqref{eq:Xib}. We give them here separately for completeness. We have for $k=1,...,m-1$
	
	\begin{equation*}\label{eq:bFCTL_prob_blocked_first}
		\mathbb{P}(X_{1} = k,S=b) =
		p_1q_1 \sum_{l=1}^k \mathbb{P}(X_{c} = l)\mathbb{P}(Y_1 = k-l) +  q_1\mathbb{P}(X_{c}=0)\mathbb{P}(Y_{1,b} = k),
	\end{equation*} and for $i=2,\dots,g_1$ and $k=1,\dots,m-1$
	\begin{equation*}\label{eq:bFCTL_prob_blocked}
		\begin{aligned}
			\mathbb{P}(X_{i} = k, S=b) = &
			\sum_{l=1}^{k} \left\{p_iq_i\mathbb{P}(X_{i-1}=l,S=u)+q_i\mathbb{P}(X_{i-1}=l,S=b)\right\}\mathbb{P}(Y_i = k-l)\\&+ q_i\mathbb{P}(X_{i-1} = 0,S=u)\mathbb{P}(Y_{i,b} = k),
		\end{aligned}
	\end{equation*}
	which provides us with the $(m-1)g_1$ additional equations. In total, we obtain a set of $m(g_1+g_2)+(m-1)g_1$ linear equations with $m(g_1+g_2)+(m-1)g_1$ unknowns, which we can solve to find the unknown $\mathbb{P}(X_{i}=l,S=u)$, for $i=1,\dots,g_1$ and $l=0,\dots,m-1$, the unknown $\mathbb{P}(X_{i}=l,S=b)$, for $i=1,\dots,g_1$ and $l=1,\dots,m-1$, and the unknown $\mathbb{P}(X_{i}=l)$, for $i=g_1+1,\dots,g_1+g_2-1$, $i=c$, and $l=0,\dots,m-1$. Due to the complicated structure of our formulas, we do not obtain a similar, easy-to-compute Vandermonde system as for the standard FCTL queue (see~\cite{van2006delay}), but a linear solver is in general able to find the unknowns (we did not encounter any numerical issues/problems in the examples that we studied).

	There are several ways to obtain the roots of $X_d(z)$ in Equation~\eqref{eq:generalbFCTL}. Because those roots are subsequently used in solving a system of linear equations, we need to find the required roots with a sufficiently high precision, certainly if $m(g_1+g_2)+(m-1)g_1$ is large. In some cases, the roots can be found analytically, e.g. in case the number of arrivals per slot has a Poisson or geometric distribution. In other cases, the roots have to obtained numerically. There are several ways to do so. An algorithm to find roots is given in~\cite{boon2019pollaczek}, Algorithm~1, while two other methods, one based on a Fourier series representation and one based on a fixed point iteration, are described in~\cite{janssen2005analytic}.
	
	\subsection{Performance measures}\label{sec:performance_measures}
	
	Now that we have a complete characterization of $X_{g_1+g_2}(z)$, we can find the PGFs of the queue-length distribution at the end of the other slots by employing Equations~\eqref{eq:X1b} up to \eqref{eq:Xg1+g2+i}. This basically implies that we can find any type of performance measure related to the queue-length distribution. As an example we find the PGF of the queue-length distribution at the end of an arbitrary slot. We denote this PGF with $X(z)$ and obtain the following expression:
	\begin{equation*}
		X(z) = \frac{1}{c}\sum_{i=1}^c X_{i}(z) .
	\end{equation*}
	
	Another important performance measure is the delay distribution. The mean of the delay distribution, $\mathbb{E}[D]$, can easily be derived from the mean queue length at the end of an arbitrary slot by means of Little's law with a time-varying arrival rate (for a proof of Little's law in this setting see e.g.~\cite{stidham1972lambda}):
	\begin{equation*}
		\mathbb{E}[D] = \frac{X^\prime(1)}{\frac{1}{c}\sum_{i=1}^cY_{i}^\prime(1)}.
	\end{equation*}
	The PGF of the delay distribution can be derived (as is done for the FCTL queue in~\cite{van2006delay}), but such a derivation is more difficult. In the regular FCTL queue, the number of slots an arriving car has to wait is deterministic when conditioned on the number of vehicles in the queue and the time slot in which the car arrives. This is not the case for the bFCTL queue as the occurrence of blockages is random. By proper conditioning on the various blocked slots and queue lengths, one can obtain the delay distribution from the distribution of the queue length. We do not pursue this here.
	
	If we want to obtain probabilities and moments from a PGF, we need to differentiate the PGF and respectively put $z=0$ or $z=1$. In our experience, this has not proven to be a problem. However, differentiation might become prohibitive in various settings, e.g. when $m(g_1+g_2)+(m-1)g_1$ becomes large or if we want to obtain tail probabilities. There are ways to circumvent such problems. If we are pursuing probabilities and do not want to rely on differentiation, we might use the algorithm developed by Abate and Whitt in~\cite{abate1992numerical} to numerically obtain probabilities from a PGF. For obtaining moments of random variables from a PGF, an algorithm was developed in~\cite{choudhury1996numerical} which finds the first $N$ moments of a PGF numerically. Essentially, this shows that, from the PGF, we can obtain any type of quantity related to the steady-state distribution of the queue length, in the form of a numerical approximation.
	
	All formulas computed in this section have been verified by comparing the numerical results with a simulation which mimics our discrete-time queueing model. More information about this simulation is given in Appendix~\ref{a:simulation}.
	
	%Materials and Methods should be described with sufficient details to allow others to replicate and build on published results. Please note that publication of your manuscript implicates that you must make all materials, data, computer code, and protocols associated with the publication available to readers. Please disclose at the submission stage any restrictions     on the availability of materials or information. New methods and protocols should be described in detail while well-established methods can be briefly described and appropriately cited.
	
	%Research manuscripts reporting large datasets that are deposited in a publicly available database should specify where the data have been deposited and provide the relevant accession numbers. If the accession numbers have not yet been obtained at the time of submission, please state that they will be provided during review. They must be provided prior to publication.
	
	%Interventionary studies involving animals or humans, and other studies require ethical approval must list the authority that provided approval and the corresponding ethical approval code.

	\section{Examples}\label{sec:results}
	
	We start in Subsection~\ref{subsec:special} with several special cases of the bFCTL queue for which we provide explicit expressions for the PGF of the overflow queue and relate those special cases to the existing literature. Subsequently, we make a comparison between the capacity obtained in the HCM~\cite{HCM} and the capacity in our model in Subsection~\ref{subsec:capacity2}. After that, we investigate the influence of several parameters on the performance measures in numerical examples. We consider performance measures like the mean and variance of the steady-state queue-length distribution, both at specific moments and at the end of an arbitrary slot, the mean delay, and several interesting queue-length probabilities. We study the influence of the $p_i$ and $q_i$ in Subsection~\ref{subsec:parameter}. In Subsection~\ref{subsec:layout}, we compare the case of turning and straight-going traffic on a single lane, as present in the bFCTL queue where blockages of all vehicles might occur, and cases where we have dedicated lanes for the right-turning and straight-going traffic where only turning vehicles are blocked. Note that we will consider each lane \emph{separately} in those examples, so there is no conflict with e.g. Remark~\ref{rem:pi}. %Afterwards, we investigate the bFCTL queue with multiple lanes without any blockages, so we study a direct extension of the regular FCTL queue to a model with multiple lanes. We do this in Subsection~\ref{subsec:multiFCTL}.
	
	\subsection{Special cases of the bFCTL queue}\label{subsec:special}
	
	We study several special cases of the bFCTL queue, e.g. cases where the bFCTL queue reduces to the FCTL queue.
	
	If $q_i=1$, an explicit expression for the PGF of the distribution of the overflow queue, $X_{g_1+g_2}(z)$, can be written down relatively easily. When it is further assumed, for the ease of exposition, that all $p_i=p$, $Y_i\overset{d}= Y$, $Y_{i,b}\overset{d} = Y_b$ and $m=1$, the following expression for $X_{g_1+g_2}(z)$ is obtained:
	\begin{equation}\label{eq:explicit_Xg_q=0}
		X_{g_1+g_2}(z) = \frac{X_n(z)}{X_d(z)},
	\end{equation} with
	\begin{align}
		\label{eq:Xn_q=0}
		X_n&(z) =z^{g_1+g_2}\sum_{i=0}^{g_2-1}\left(\frac{Y(z)}{z}\right)^{g_2-i-1}\left(1-\frac{Y(z)}{z}\right)\mathbb{P}(X_{g_1+i}=0)+ \nonumber\\
		& z^{g_1}Y(z)^{g_2}\sum_{i=0}^{g_1-1}\nonumber \Bigg\{ \mathbb{P}(X_{i}=0,S=u)\Bigg[\left(Y_b(0) - (1-p)\frac{Y(z)}{z}\right)\left((1-p)\frac{Y(z)}{z}\right)^{g_1-i-1}+\\&\left(Y_b(z)-Y_b(0)-p Y(z)\right)Y(z)^{g_1-i-1}\Bigg]+\nonumber\\
		&p Y(z)^{g_1-i}\sum_{j=0}^{i-1} \mathbb{P}(X_{j}=0,S=u)\left(Y_b(0) - (1-p)\frac{Y(z)}{z}\right)\left((1-p)\frac{Y(z)}{z}\right)^{i-j-1}\Bigg\},
	\end{align}
	where $\mathbb{P}(X_{0}=0,S=u)$ is to be interpreted as $\mathbb{P}(X_{c}=0)$, and
	\begin{equation}\label{eq:Xd_q=0}
		X_d(z) = z^{g_1+g_2} - \left(\left(1-p\right)^{g_1}+p z^{g_1}\sum_{i=0}^{g_1-1}\left(\frac{1-p}{z}\right)^i\right) Y(z)^c.
	\end{equation}
	The reason that we provide an explicit formula for this particular case is that this formula is significantly easier than the formula in the case where $q_i< 1$ for one or more $i=1,\dots,g_1$. The stability condition (cf. Algorithm~\ref{alg:stability} in Subsection~\ref{subsec:capacity}) for this example is relatively easy to derive and reads as follows:
	\begin{equation*}
		\begin{cases}
			\mu c < g_1+g_2, & \textrm{if } p = 0,\\
			\mu c < g_2, & \textrm{if } p = 1,\\
			\mu c < g_2+\left(1-(1-p)^{g_1}\right)\frac{1-p}{p}, & \textrm{otherwise},
		\end{cases}
	\end{equation*}
	where $\mu$ is the mean arrival rate per slot, i.e. $\mu=\mathbb{E}[Y]$.
	This can be understood as follows: if $p=0$ there are no turning vehicles and we obtain the regular FCTL queue with green period $g_1+g_2$. If $p=1$ all vehicles are turning vehicles and there are no departures during the first part of the green period because $q_i=1$, so we obtain the FCTL queue with green period $g_2$. The other case can be understood as follows: on the left-hand side we have the average number of arrivals per cycle whereas on the right-hand side we have the average number of slots available for delayed vehicles to depart. Indeed, on the right-hand side we have $g_2$, the number of green slots during the second part of the green period which are all available for vehicles to depart, and the number of green slots available for departures during the first green period:
	\[
	\sum_{i=1}^{g_1}(1-p)^i = \left(1-(1-p)^{g_1}\right)\frac{1-p}{p}.
	\]
	
	If $p_i=0$ for all $i$, i.e. there are no blockages occurring at all (regardless of the $q_i$), the FCTL queue with multiple lanes (with green period $g=g_1+g_2$) is obtained. Note that we do not have to include the state $S$, because there are no blockages of batches of vehicles. If $m=1$, we obtain the regular FCTL queue as studied in e.g.~\cite{van2006delay}. This can e.g. be observed when putting $p_i=0$ and $m=1$ in Equations~\eqref{eq:explicit_Xg_q=0}, \eqref{eq:Xn_q=0}, and \eqref{eq:Xd_q=0}. The expression for $X_{g_1+g_2}(z)$ or, alternatively, $X_g(z)$ is (after rewriting):
	\begin{align}\label{eq:FCTL}
		X_g(z) & = \frac{(z-Y(z))z^{g-1}\sum_{i=0}^{g-1}\mathbb{P}(X_i=0)\left(\frac{Y(z)}{z}\right)^{g-i-1}}{z^g-Y(z)^c},
	\end{align}
	where $\mathbb{P}(X_{0}=0)$ is to be interpreted as $\mathbb{P}(X_{c}=0)$. For general $m$, we have the following formula:
	\begin{equation}
		X_g(z) = \frac{z^{mg}\sum_{i=0}^{g-1}\sum_{l=0}^{m-1} \mathbb{P}(X_{i}=l)\left(1-\frac{Y(z)}{z^{m-l}}\right)\left(\frac{Y(z)}{z}\right)^{g-i-1}}{z^{mg}-Y(z)^c},
	\end{equation}
	where the $\mathbb{P}(X_{0}=l)$, $l=0,\dots,m-1$, are to be interpreted as $\mathbb{P}(X_{c}=l)$. The stability condition for this case can be verified to be
	\[
	\mu c < m g
	\]
	which is in accordance with Algorithm~\ref{alg:stability}.
	
	It can also be verified that the bFCTL queue reduces to the regular FCTL queue with green time $g = g_2$ and red time $r+g_1$, if $p_i = 1$ and $q_i = 1$.
	
	We note that for the FCTL queue with a single lane and no blockages (i.e. $p_i=0$ or $p_i=1$ and $q_i=1$) there is an alternative characterization of the PGF in terms of a complex contour integral, see~\cite{boon2019pollaczek}. It remains an open question whether such a contour-integral representation exists for the bFCTL with multiple lanes, as the polynomial structure in terms of $Y(z)/z$ as present in Equation~\eqref{eq:FCTL} is not present in the general bFCTL queue. This feature of the FCTL queue seems essential to obtain a contour-integral expression as is done in~\cite{boon2019pollaczek}.
	
	In \cite{boon2019pollaczek}, a decomposition result is presented in Theorem~2. It shows that several related queueing processes can in fact be decomposed in the independent sum of the FCTL queue and some other queueing process. It is likely that the bFCTL queue with multiple lanes allows for some of those generalizations as well. We mention randomness in the green and red time distributions as a relevant potential extension.
	
	\subsection{Capacity} \label{subsec:capacity2}
	
	In order to compare our model and the existing literature (focusing on the HCM~\cite{HCM}), we provide several examples in this subsection.
	
	The formula for the capacity of a permitted right-turn lane in a shared lane in the HCM is
	\begin{equation*}
		s_{sr} = \frac{s_{th}}{1+P_r\left(\frac{E_R}{f_{Rpb}}-1\right)},
	\end{equation*}
	cf.~\cite{HCM} equation (31-105). Here, $s_{sr}$ is the saturation flow of the shared lane, $s_{th}$ the saturation flow of an exclusive through lane, $P_r$ the right-turning portion of vehicles, $E_R$ the equivalent number of through vehicles for a protected right-turn vehicle and $f_{Rpb}$ is the bicycle-pedestrian adjustment factor for right-turn groups. The latter is defined as the average amount of time during the green period during which right-turning vehicles are not blocked, i.e., in our model, there are no pedestrians crossing. There is a procedure provided in the HCM to compute this factor, but in our model this simply corresponds to the $q_i$ and we will determine the $f_{Rpb}$ factor on the $q_i$. Further, in order to make a comparison with our model, we turn the saturation flow of the shared lane into a number of vehicles per cycle.
	More concretely, we choose the green period to be $30$ seconds, split into the two phases as follows: $g_1=20$ and $g_2=10$. We pick the cycle length to be $90$ seconds, the time slots to have length $2$ seconds and we focus on a single shared lane, so we have at most $1$ vehicle departing per time slot. Further, we choose the right-turning portion vehicles to be $1$ or $0.9$ in our examples. Lastly, for the HCM formula, we assume that vehicles heading straight have a crossing time of $1$ second. To account for this effect in our bFCTL model, we use the correction discussed in Remark~\ref{rem:capacity}. In this example we have:
	\[
	m^* = p_i m_\textit{turn} + (1-p_i)m_\textit{through} = p_i\times 1 + (1-p_i)\times 2.
	\]
	This enables us to compute the capacity in our model and in the HCM up to the $q_i$.
	
\begin{figure}[h!]
		\centering
		\begin{tabular}{cc}
			\includegraphics[width=0.45\textwidth]{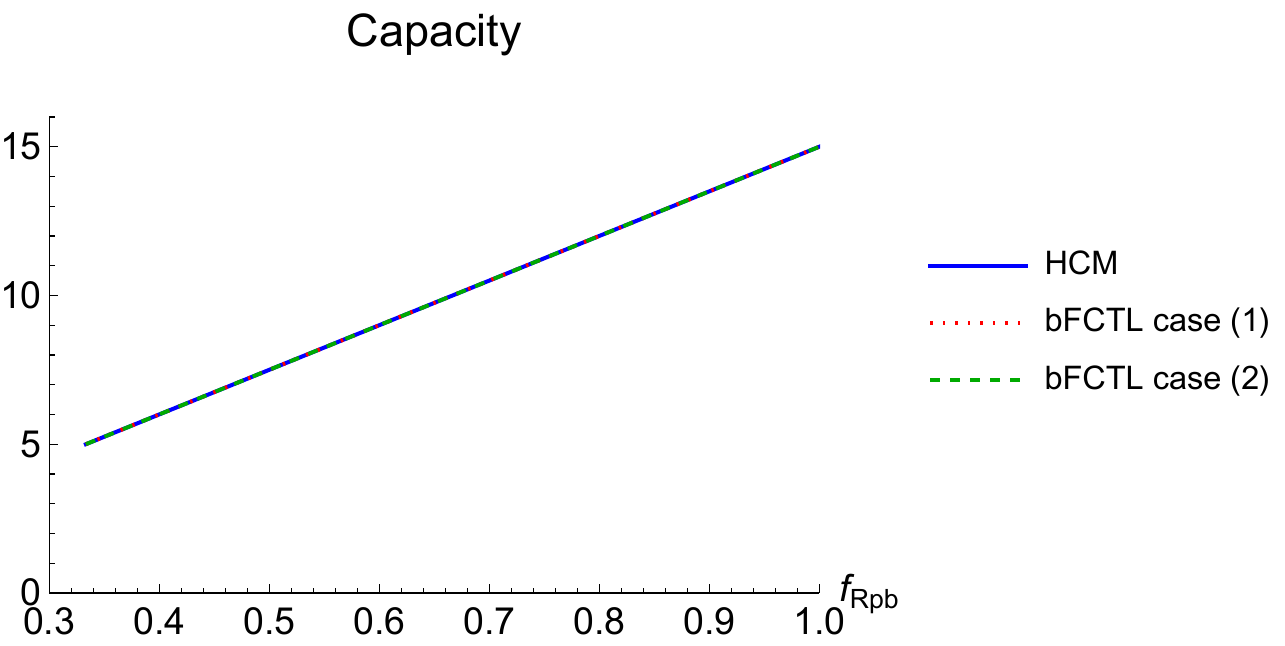}
			&
			\includegraphics[width=0.45\textwidth]{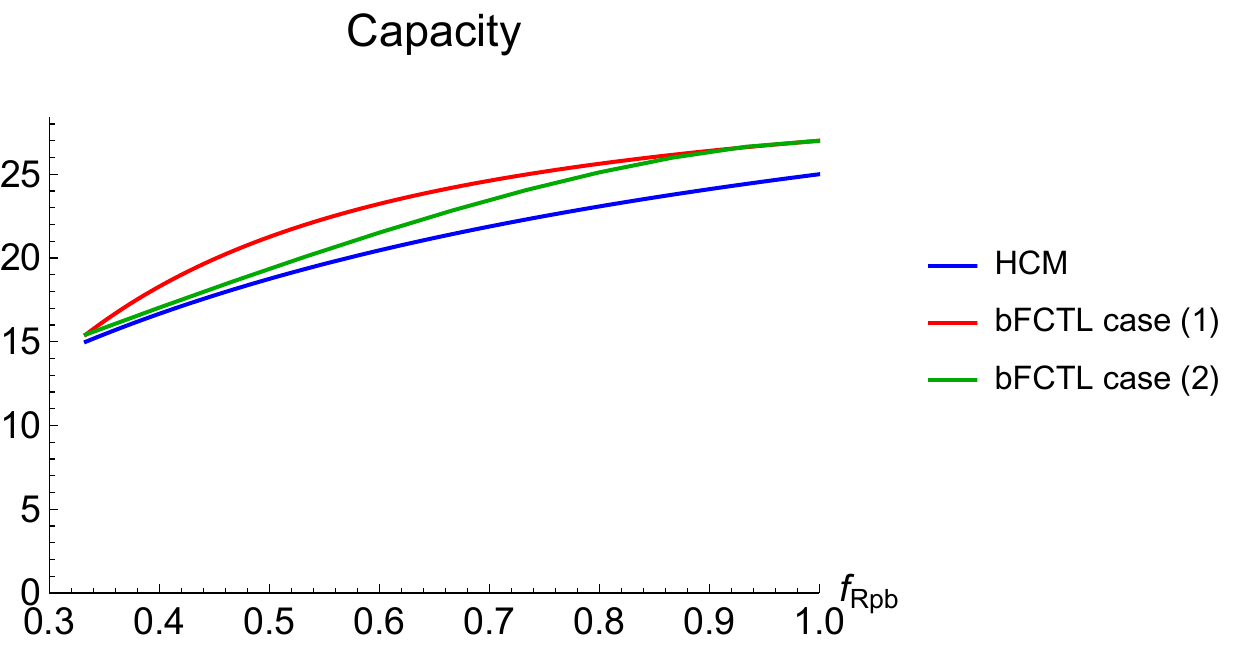}\\
			(a) & (b) \\[1ex]
		\end{tabular}
		\caption{Capacity in vehicles per cycle for the example according to the HCM and to the bFCTL queue with two different choices for the $q_i$ (case $(1)$ and case $(2)$) as detailed in the text. We have $p_i=1$ in (a) and $p_i=0.2$ in (b).}
		\label{f:cap1}
	\end{figure}

	We first focus on the cases with $p_i=1$ and we display the capacity according to the HCM in Figure~\ref{f:cap1}(a). Note that $f_{Rpb}$ is at least $1/3$ because $g_1=20$ and $g_2=10$, implying that during at least a part $1/3$ of the cycle, turning vehicles are not blocked. In Figure~\ref{f:cap1}(a) we also depict two capacities according to the bFCTL queue. In case $(1)$ we assume that all the $q_i$ are the same and are chosen in such a way that the $f_{Rpb}$ in the HCM formula is matched. E.g. in case $f_{Rpb}=1/3$, we choose $q_i=0$ as there are no pedestrians and in case $f_{Rpb}=2/3$, we choose $q_i=1/2$. In case $(2)$, we consider a step function for the $q_i$ such that
	\begin{equation*}
		q_i = \begin{cases}
			1 & \textrm{if } i < k\\
			0 &  \textrm{if } i > k\\
			k^* & \textrm{otherwise,}
		\end{cases}
	\end{equation*}
	for some values of $k$ and $k^*$ such that the $q_i$ match with the value for $f_{Rpb}$ that is used in the formula for the HCM.
	
	Figure~\ref{f:cap1}(a) makes sense: if $f_{Rpb}$ is for example equal to $1$, there are no pedestrians crossing (i.e. $q_i=0$), and then the number of vehicles departing per cycle is $(g_1+g_2)/2=15$. The capacity according to the bFCTL queue when $p_i=1$ is equal to (after simplification)
	\begin{equation}\label{eq:cap}
		\frac{g_1+g_2}{2} - \sum_{i=1}^{g_1} q_i.
	\end{equation}
	\textbf{This shows that when $\sum_{i=1}^{g_1} q_i$ is translated into the factor $f_{Rpb}$ in the HCM, we have an identical capacity.} E.g. if the $q_i=0$, then also in the bFCTL queue, the capacity is equal to $15$ vehicles per cycle. Equation~\eqref{eq:cap} also indicates that it does not matter in which slots the pedestrians are crossing if $p_i=1$ (when looking at the capacity). In this case, the $q_i$ only influence the capacity through their sum, however in general the individual $p_i$ and $q_i$ have an impact on the capacity (and the queue-length process). Similar observations hold if $p_i=0$, i.e. there are no turning vehicles.
	
	\textbf{If the $p_i$ are not equal to $1$, there are differences between the capacity in the HCM and the bFCTL queue.} We study an example where $p_i=0.2$. The results are depicted in Figure~\ref{f:cap1}(b). The values for the capacity obtained with the function in the HCM are slightly lower than the values that we obtain in both cases of the bFCTL queue.
	
	In contrast with the previous example, there are differences between all three choices which relate to various causes. %One explanation for the difference between the HCM and the bFCTL queue is that there is a difference in the crossing times between right-turning vehicles and vehicles heading straight in the HCM, which is not captured in the bFCTL queue. We note that this difference is more profound for lower values of $p_i$ (based on more experiments). If we would adjust for this in the HCM (e.g. by adjusting $E_R$ and the $s_{th}$), our results still differ from the HCM due to the more realistic description in the bFCTL queue, but the absolute differences decrease. Vice versa, we could adjust the framework of bFCTL queue and abandon the direct relation between slots and vehicle departures. In more detail, departing vehicles could take an integer number of slots for the duration of a departure (e.g. either $1$ or $2$ slots), as is demonstrated for the FCTL model in Section 3.3 in~\cite{maesnetworks}. In this section,  the state of the Markov chain is extended with an extra dimension,  represented by an additional variable that keeps track of whether a vehicle is actually allowed to depart in a given time slot (the departure variable in~\cite{maesnetworks}). Although  the  idea seems  relatively  straightforward,  a  realistic  implementation  would  require making significant changes to the Markovian bFCTL model and is considered to be beyond the scope of this paper. Another benefit of such an approach would be that it probably enables the modelling of different saturation flows among the lanes, which would further add to the practical relevance of our approach.
	The main reason for the occurring difference between cases $(1)$ and $(2)$ in the bFCTL queue, is that the individual $q_i$ are determining the capacity rather than the total value of the $q_i$'s alone as was the case when $p_i=1$. \textbf{Here we thus see that our detailed description of the queueing model in terms of slots is necessary to fully understand the capacity (and, more generally, the queueing process).}
	
	In this subsection we have been working under several assumptions. If one would, e.g., also incorporate start-up delays as is done in~\cite{shaoluen2020random}, we would see that the capacity in the HCM results in an overestimation of the capacity as is more generally observed~\cite{shaoluen2020random}. We also expect that the distribution of the $q_i$ over the different slots has a bigger impact on the capacity and queueing process if start-up delays are incorporated. Implementing such effects into our model is possible (probably in a similar way as including a departure variable as discussed above), but is beyond the scope of the present paper.
	
	\subsection{The bFCTL queue with turning vehicles and pedestrians}\label{subsec:parameter}
	
	In this subsection, we study the bFCTL queue with a single lane, so $m=1$. The setting in this subsection is as depicted in Figure~\ref{fig:vis}(b). We mainly focus on the distribution of $X_{g_1+g_2}$, to which we refer as the overflow queue, as this is the distribution from which some interesting performance measures can be derived.
	This distribution reflects the probability distribution of the queue size at the moment that the green light switches to a red light. We also briefly consider some other performance measures.
	
	\subsubsection{Influence of the number of turning vehicles}\label{sec:influenceofturning}
	
\begin{figure}[h!]
		\centering
		\begin{tabular}{cc}
			\includegraphics[height=4.7cm]{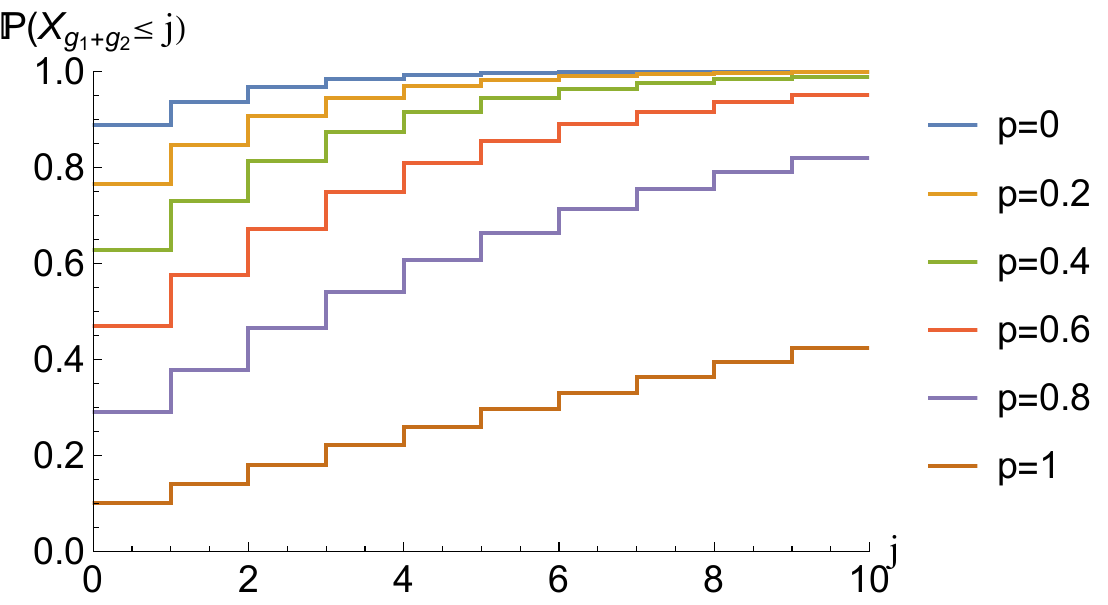}
			&
			\includegraphics[height=4.7cm]{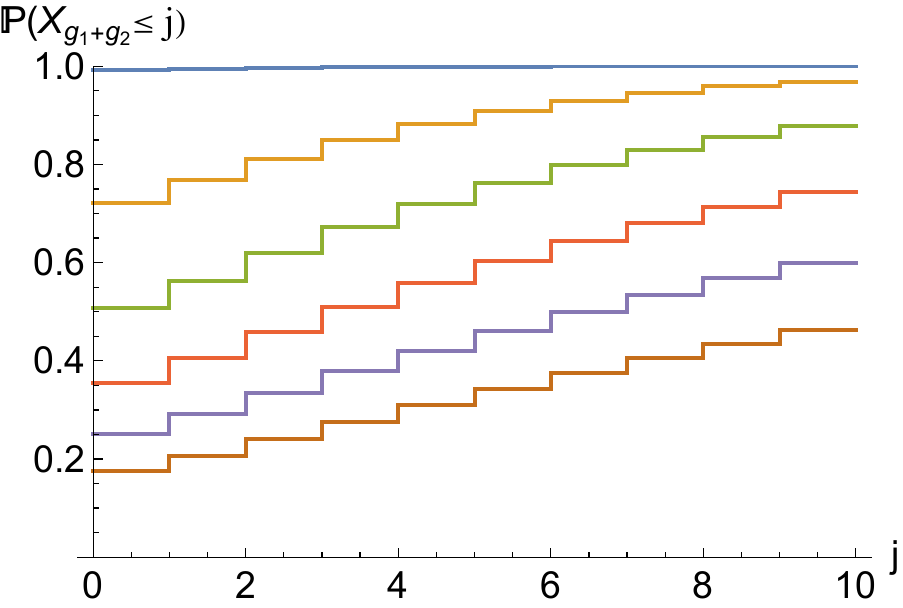}\\
			(a) & (b) \\[1ex]
		\end{tabular}
		\caption{Cumulative Distribution Function (CDF) of the overflow queue for various values of $p_i=p$, $q_i=q=1$, and Poisson arrivals with mean $0.39$. In (a) we have $g_2=r=2g_1=4$ and in (b) we have $g_2=r=2g_1=20$.}
		\label{fig:ex3}
	\end{figure}

	First, we vary the fraction of right-turning vehicles $p_i$ and study its influence on $X_{g_1+g_2}$. We choose the $p_i$ to be the same for each $i$, so we have $p_i=p$, and we vary $p$. We choose the value of the $q_i=q$ to be $1$, so there are always pedestrians on the pedestrian crossing during the first part of the green period with length $g_1$. In this way, we can effectuate the influence of the fraction of turning vehicles on the performance measures. Further, we choose $g_1$ to be either $2$ or $10$ and we choose $g_2=r=2g_1$. The arrival process is taken to be Poisson with mean $0.39$. Note that  the lane is close to its point of saturation, because the capacity can be shown to be equal to $0.4$.
	We display results for $\mathbb{P}(X_{g_1+g_2}\leq j)$ for $j=0,\dots,10$ in Figure~\ref{fig:ex3}.

	As can be observed from Figure~\ref{fig:ex3}, \textbf{the fraction of turning vehicles may dramatically influence the number of queueing vehicles}. There is virtually no queue at the end of the green period when there are no turning vehicles ($p=0$), whereas in \textbf{more than 50\%} of the cases there is a queue of at least $10$ vehicles at the end of the green period when all vehicles are turning vehicles ($p=1$). The blockages of the turning vehicles in the latter case effectively reduce the green period by a factor $1/3$ in our examples (as $q=1$), which causes the huge difference in performance. We note that the distribution of $X_{g_1+g_2}$ coincides with the overflow queue distribution in the FCTL queue when $p=0$ (when we take $g_1+g_2$ as the green period and $r$ as the red period in the FCTL queue) and when $p=1$ and $q=1$ (with $g_2$ the green period and $r+g_1$ the red period).
	
	When comparing Figures~\ref{fig:ex3}(a) and ~\ref{fig:ex3}(b), we see that \textbf{the influence of $p$ is not uniform across the two examples}. In case $p=0$ or $p=1$, the probability of a large overflow queue is larger for the case where $g_1=2$. This might be clarified by noting that a larger cycle reduces the amount of within-cycle variance which reduces the probabilities of a large queue length. If $0<p<1$ this does not seem to be the case. This might be due to the fact that a relatively big part of the first green period is eaten away by turning vehicles that are blocked when $g_1=10$. For example, when $p>0$ and the first vehicle is a turning vehicle, immediately the entire period $g_1$ is wasted because $q=1$. This is of course also the case when $g_1=2$, but the blockage is resolved sooner and during the second part of the green period the blocked vehicle may depart relatively soon in comparison with the case where $g_1=10$. %It might be the case that the within-cycle variance is bigger in the case $g_1=10$ and $0<p<1$ than in the case $g_1=2$ and $0<p<1$.
	
	\begin{figure}[h!]
		\centering
		\begin{tabular}{cc}
			\includegraphics[width=0.47\textwidth]{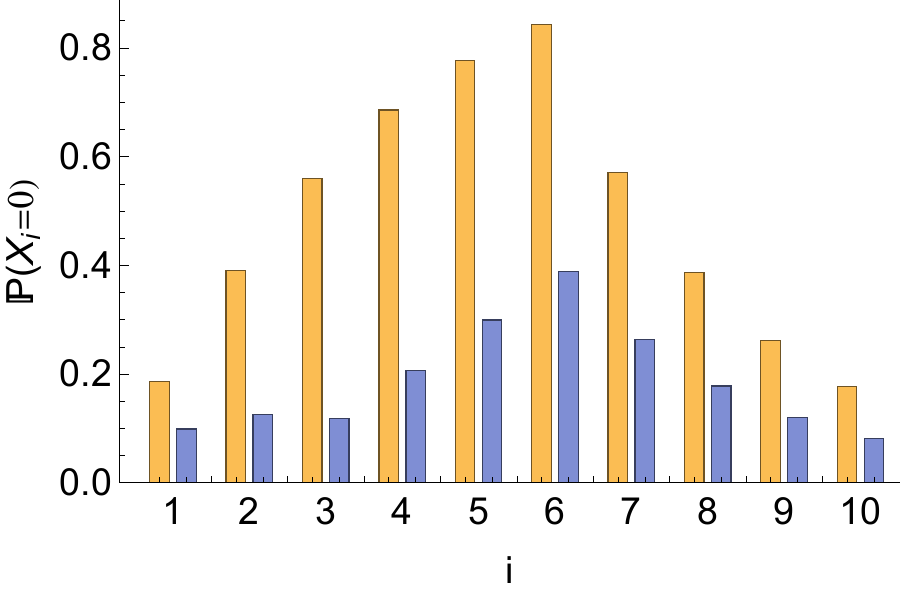}
			&
			\includegraphics[width=0.47\textwidth]{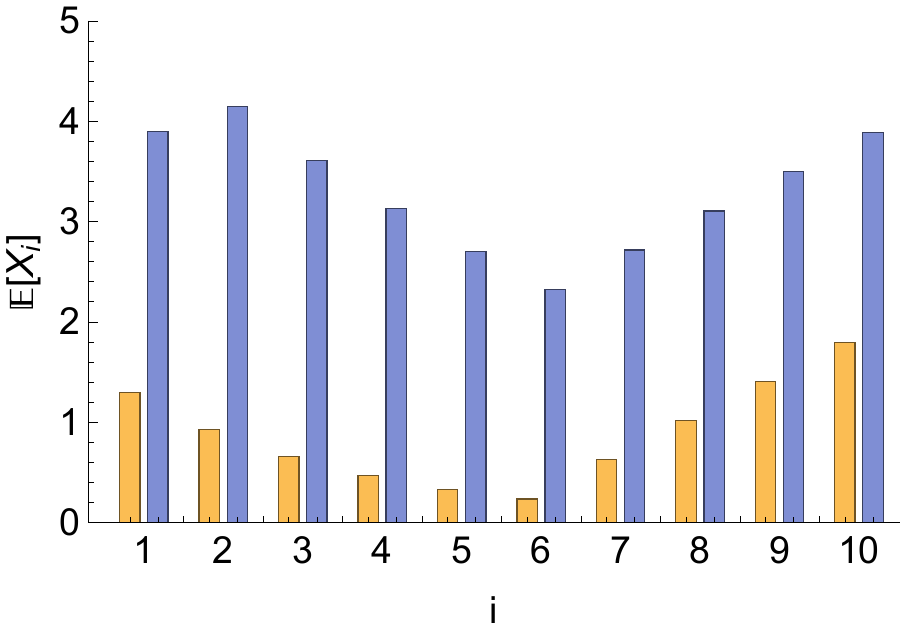}\\
			(a) & (b) \\[1ex]
		\end{tabular}
		\caption{In~(a) $\mathbb{P}(X_{i}=0)$ for slot number $i=1,\dots,10$ is displayed for two different values of $p_i$, where orange corresponds to $p_i=p=0$ and blue to $p_i=p=0.6$, with $2g_1=g_2=r = 4$, $q_i=q=1$, and with Poisson arrivals with mean $0.39$. In~(b) the same two examples are studied, but the mean queue length $\mathbb{E}[X_{i}]$ at the end of slot $i$ is shown. }
		\label{fig:ex3xk}
	\end{figure}
	
	In Figure~\ref{fig:ex3xk}(a), we see the probability of an empty queue after slot $i$, where $i=1,2,\dots,c$, for two different values of $p$. For the case $p=0$ (in orange) we have a monotone increasing sequence of probabilities during the green period as one would expect: this setup corresponds to a regular FCTL queue and once the queue empties during the green period, it stays empty. We see that for the case $p=0.6$ (in blue) the probabilities of an empty queue after slot $i$ are much lower (as there are more turning vehicles which might be blocked and hence cause the queue to be non-empty). In fact, the probability of an empty queue even decreases when going from slot $2$ to slot $3$. This can be clarified by the fact that the queue might start building again even when the queue is (almost) empty: e.g. if the queue is empty during the first green period and there is an arrival of a turning vehicle, that vehicle will be blocked as $q=1$ in which case the queue is no longer empty.
	
	The same type of behaviour is reflected in the mean queue length at the end of a slot, as can be observed in Figure~\ref{fig:ex3xk}(b). Even though the green period already started, the queue in the example with $p=0.6$ still grows (in expected value) during the first part of the green period, see the first two blue bars. This is caused by the fact that vehicles might be blocked, \textbf{which demonstrates the possibly severe impact of blocked vehicles on the performance of the system}.
	
	\subsubsection{Influence of the pedestrians}
	
\begin{figure}[h!]
		\centering
		\begin{tabular}{cc}
			\includegraphics[height=4.5cm]{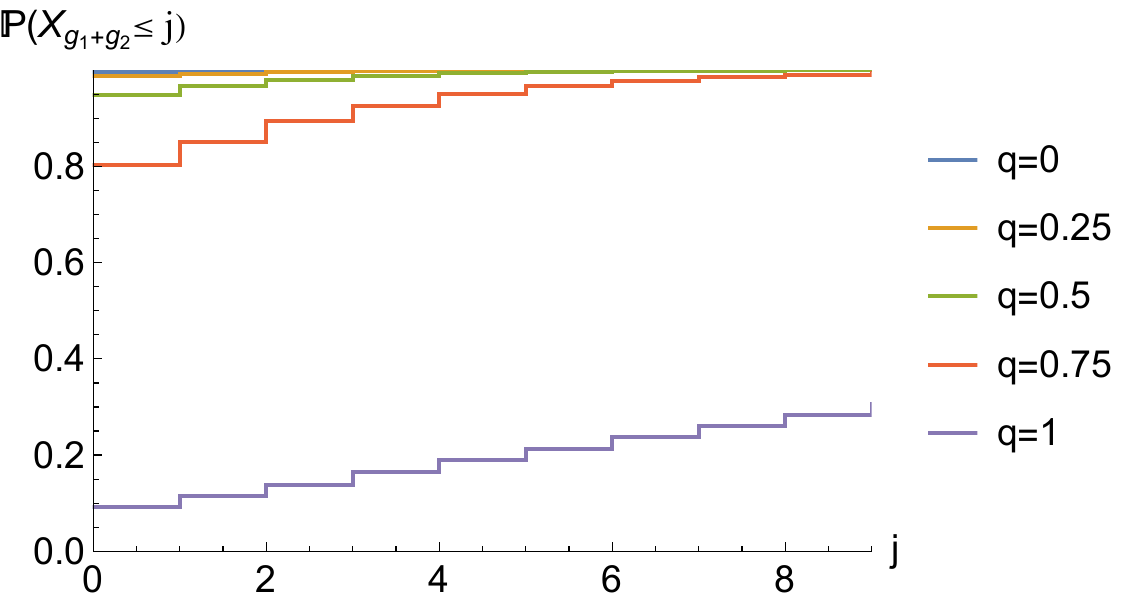}
			&
			\includegraphics[height=4.5cm]{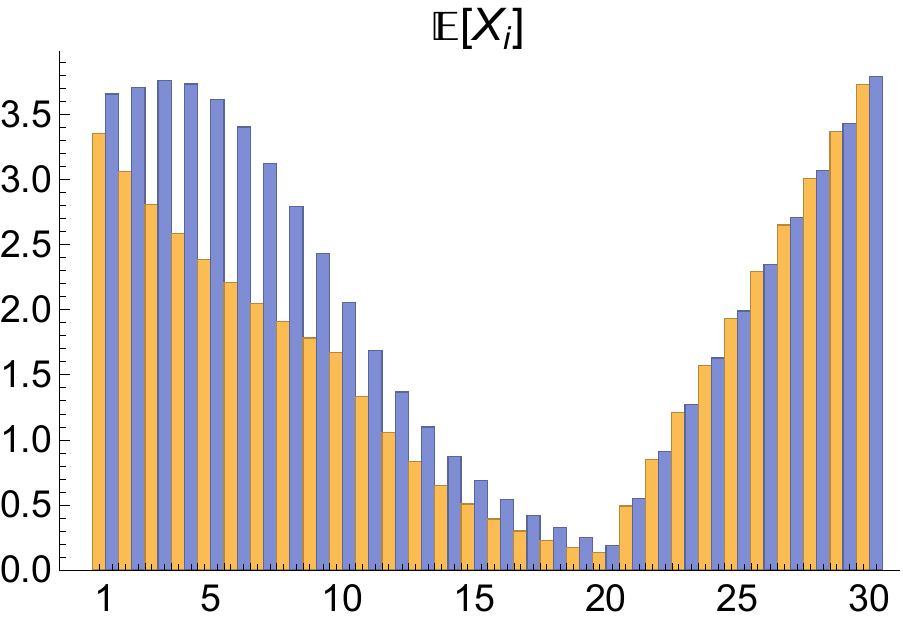}\\
			(a) & (b) \\[1ex]
		\end{tabular}
		\caption{In (a) the CDF of the overflow queue is displayed for various values of the $q_i$ with all $q_i=q$ the same, $p_i=p=0.5$, Poisson arrivals with mean $0.36$, and $g_1=g_2=r=10$. In (b) the $\mathbb{E}[X_{i}]$ are compared for slot number $i=1,\dots,30$ with in orange $q_i=0.5$ and in blue $q_i=1-(i-1)/g_1$ for $i=1,\dots,g_1$. Further, it is assumed that $p_i=p=0.5$, that the number of arrivals in each slot follows a Poisson distribution with mean $0.36$, and that $g_1=g_2=r=10$.}
		\label{fig:ex4xk}
	\end{figure}

	Secondly, we investigate the influence of the presence of pedestrians by studying various values for the $q_i$. A high value of the $q_i$ corresponds to a high density of pedestrians as $q_i$ corresponds to the probability that a turning vehicle is not allowed to depart during the first green period. Conversely, a low value of the $q_i$ corresponds to a low density of pedestrians and a relatively high probability of a turning vehicle departing during the first green period. We choose $p_i=p=0.5$ and take $g_1=g_2=r=10$. We take Poisson arrivals with mean $0.36$. We study one set of examples where the $q_i$ are constant over the various slots, see Figure~\ref{fig:ex4xk}(a). We also study the influence of the dependence of the $q_i$ on $i$ by investigating two cases with all parameters as before in Figure~\ref{fig:ex4xk}(b). In one case we take $q_i=0.5$ for all $i$, but in the other case we take $q_i=1-(i-1)/g_1$. The latter case reflects a decreasing number of pedestrians blocking the turning flow of vehicles during the first part of the green period.
	
	We note that it is \textbf{important to estimate the correct blocking probabilities $q_i$ from data}, when applying our analysis to a real-life situation \textbf{as the $q_i$ have an impact on the performance measures}. In Figure~\ref{fig:ex4xk}(a), we clearly see that the more pedestrians, the longer the queue length at the end of the green period is. Indeed, if there are more pedestrians, there are relatively many blockages of vehicles which causes the queue to be relatively large.
	
	Moreover, \textbf{it is important to capture the dependence of the $q_i$ on the slot $i$ in the right way}, see Figure~\ref{fig:ex4xk}(b). Even though, on average over all slots, the mean number of pedestrians present is similar in the two cases, we see a clear difference between the two examples. In the case with decreasing $q_i$ (in blue), we see an initial increase of the mean queue length during the first green slots of the cycle, caused by a relatively large fraction of turning vehicles ($p=0.5$) \emph{and} a high value of $q_i$. This is not the case in the other example where $q_i=0.5$ for all $i$. After some slots of the first green period, the decrease in the mean queue length is quicker for the example where the $q_i$ decrease when $i$ increases, which can (at least partly) be explained by the decreasing $q_i$. During the remaining part of the cycle, the queue in front of the traffic light behaves more or less the same in both examples and even the mean overflow queue, $\mathbb{E}[X_{g_1+g_2}]$, is not that much different for the two examples. This implies, as can also be observed in Figure~\ref{fig:ex4xk}(b), that the mean queue length during the red period is comparable as well for our setting. This does not hold for the mean queue length at the end of an arbitrary slot and the mean delay, because of the differences in the queue length during the first part of the green period.
	
	\subsection{Shared right-turn lanes and dedicated lanes}\label{subsec:layout}
	
\begin{figure}[h!]
		\centering
		\begin{tabular}{ccc}
			
			\begin{tikzpicture}[scale=0.1,rotate=90]
				\draw[black,fill=lightgray](30,20) rectangle (40,50);
				\draw[black,fill=lightgray](0,30) rectangle (50,40);
				\draw[thick,white,dash pattern=on 7 off 4](0,35) to (50,35);
				%\draw[thick,white](0,40) to (50,40);
				\draw[thick,white,dash pattern=on 7 off 4](35,0) to (35,60);
				\draw[lightgray,fill=lightgray](30,30) rectangle (40,40);
				\draw[white,->,thick](24,32.5) to [out = 0, in = 120] (28,30.5);
				\draw[white,->,thick](24,32.5) to [out = 0, in = 180] (28.5,32.5);
				%\draw[green,fill=green](30,30) rectangle (30.5,35);
				\draw[white,fill=white](31,25.5) rectangle (32,28.5);
				\draw[white,fill=white](33,25.5) rectangle (34,28.5);
				\draw[white,fill=white](36,25.5) rectangle (37,28.5);
				\draw[white,fill=white](38,25.5) rectangle (39,28.5);
			\end{tikzpicture}
			&\hspace{-1cm} \begin{tikzpicture}[scale=0.1,rotate=90]
				\draw[black,fill=lightgray](30,20) rectangle (40,60);
				\draw[black,fill=lightgray](0,30) rectangle (60,50);
				\draw[thick,white,dash pattern=on 7 off 4](0,35) to (60,35);
				\draw[thick,white,dash pattern=on 7 off 4](0,45) to (60,45);
				\draw[thick,white](0,40) to (60,40);
				\draw[thick,white,dash pattern=on 7 off 4](35,0) to (35,60);
				\draw[lightgray,fill=lightgray](30,30) rectangle (40,50);
				\draw[white,->,thick](24,32.5) to [out = 0, in = 120] (28,30.5);
				\draw[white,->,thick](24,37.5) to [out = 0, in = 180] (28.5,37.5);
				%\draw[green,fill=green](30,30) rectangle (30.5,40);
				\draw[white,fill=white](31,25.5) rectangle (32,28.5);
				\draw[white,fill=white](33,25.5) rectangle (34,28.5);
				\draw[white,fill=white](36,25.5) rectangle (37,28.5);
				\draw[white,fill=white](38,25.5) rectangle (39,28.5);
				
				%						\draw[white,->,thick](24,32.5) to [out = 0, in = 120] (28,30.5);
				%			\draw[white,->,thick](24,32.5) to [out = 0, in = 180] (28.5,32.5);
				%			\draw[white,->,thick](24,37.5) to [out = 0, in = 240] (28,39.5);
				
			\end{tikzpicture}
			&\hspace{-1cm} \begin{tikzpicture}[scale=0.1,rotate=90]
				\draw[black,fill=lightgray](30,20) rectangle (40,60);
				\draw[black,fill=lightgray](0,30) rectangle (60,50);
				\draw[thick,white,dash pattern=on 7 off 4](0,35) to (60,35);
				\draw[thick,white,dash pattern=on 7 off 4](0,45) to (60,45);
				\draw[thick,white](0,40) to (60,40);
				\draw[thick,white,dash pattern=on 7 off 4](35,0) to (35,60);
				\draw[lightgray,fill=lightgray](30,30) rectangle (40,50);
				\draw[white,->,thick](24,32.5) to [out = 0, in = 120] (28,30.5);
				\draw[white,->,thick](24,32.5) to [out = 0, in = 180] (28.5,32.5);
				\draw[white,->,thick](24,37.5) to [out = 0, in = 180] (28.5,37.5);
				%\draw[green,fill=green](30,30) rectangle (30.5,40);
				\draw[white,fill=white](31,25.5) rectangle (32,28.5);
				\draw[white,fill=white](33,25.5) rectangle (34,28.5);
				\draw[white,fill=white](36,25.5) rectangle (37,28.5);
				\draw[white,fill=white](38,25.5) rectangle (39,28.5);
			\end{tikzpicture}
			\\
			\hspace{-0.9cm}\scriptsize (a)  %Platoon forming approaching the intersection
			&
			\hspace{-2.85cm}\scriptsize (b)  %Crossing the intersection in platoons
			&
			\hspace{-2.85cm}\scriptsize (c)  %Crossing the intersection in platoons
		\end{tabular}
		\caption{The various lane configurations considered in Subsection~\ref{subsec:layout}. In~(a) we have a single lane with a shared right-turn lane. In~(b) we have two dedicated lanes: one for straight-going vehicles and one for right-turning traffic, whereas in~(c) we have a two-lane setup with one lane for straight-going vehicles only and a shared right turn.}
		\label{fig:vis2}	
	\end{figure}
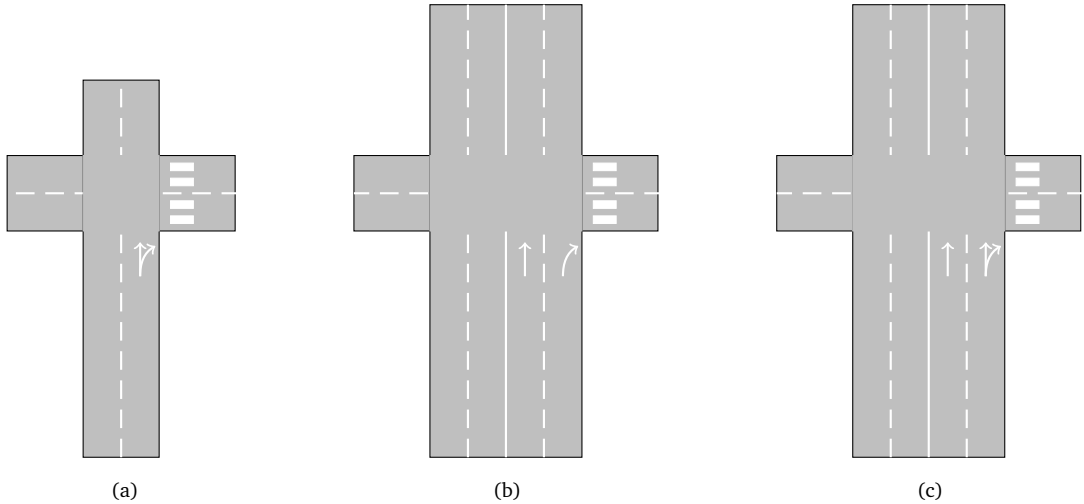

	We continue with a study of several numerical examples that focus on the differences between shared right-turn lanes and dedicated lanes for turning traffic. We do so in order to provide relevant insights in the benefit of splitting the vehicles in different streams. Firstly, we study the difference between a single shared right-turn lane (as visualized in Figure~\ref{fig:vis2}(a)) and a case where the straight-going and turning vehicles are split into two different lanes. In the latter case, we thus have two lanes, one for the straight-going traffic and one for the turning traffic (as visualized in Figure~\ref{fig:vis2}(b)) which we can analyze as two separate bFCTL queues.
	
	Secondly, we compare two two-lane settings. The first is visualized in Figure~\ref{fig:vis2}(b), while the other is a two-lane scenario where one lane is a dedicated lane for straight-going traffic and the other is a shared right-turn lane as depicted in Figure~\ref{fig:vis2}(c). We thus allow for straight-going traffic to mix with some of the right-turning vehicles in the latter case. We do so in order to make sure that the shared right-turn lane together with the lane for vehicles heading straight has the same capacity as the two lanes where the two streams of vehicles are split (as opposed to the first example in this subsection). In both two-lane scenarios we, again, analyze the two lanes as two separate bFCTL queues.
	
	\subsubsection{One lane for the shared right-turn}
	
\begin{figure}[h!]
		\centering
			\includegraphics[width=0.7\textwidth]{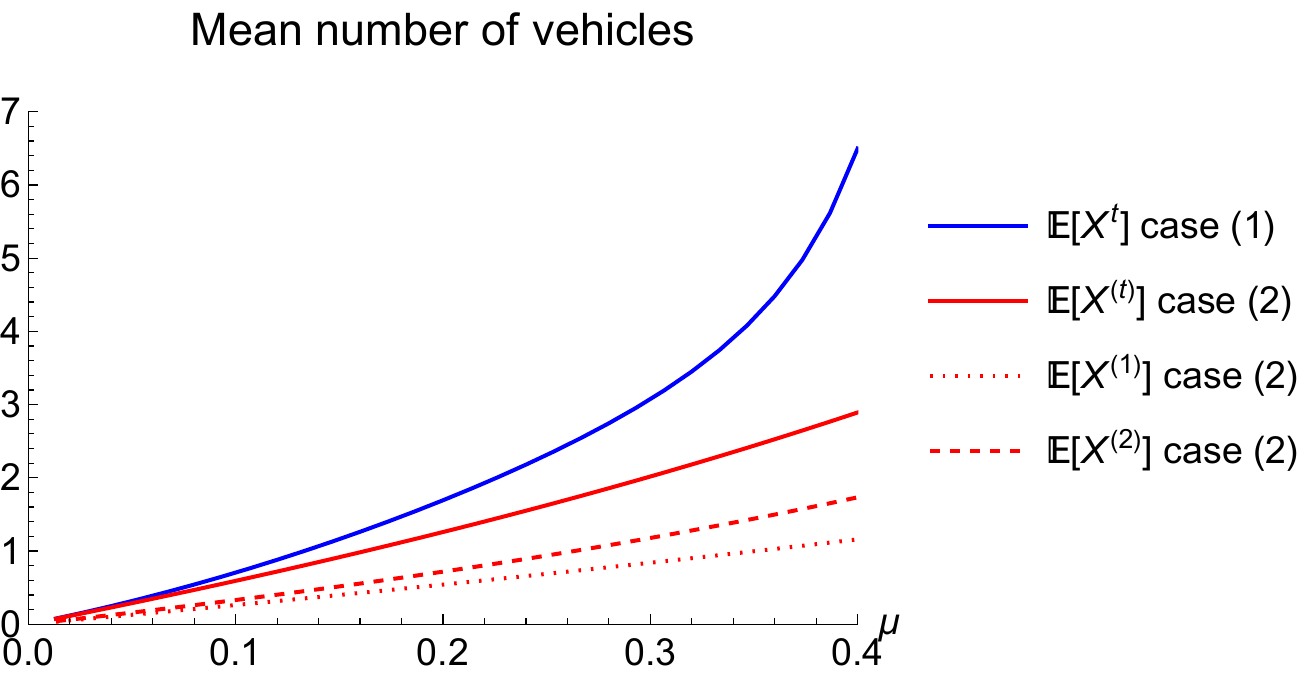}\\
			(a)\\[1ex]
			\includegraphics[width=0.7\textwidth]{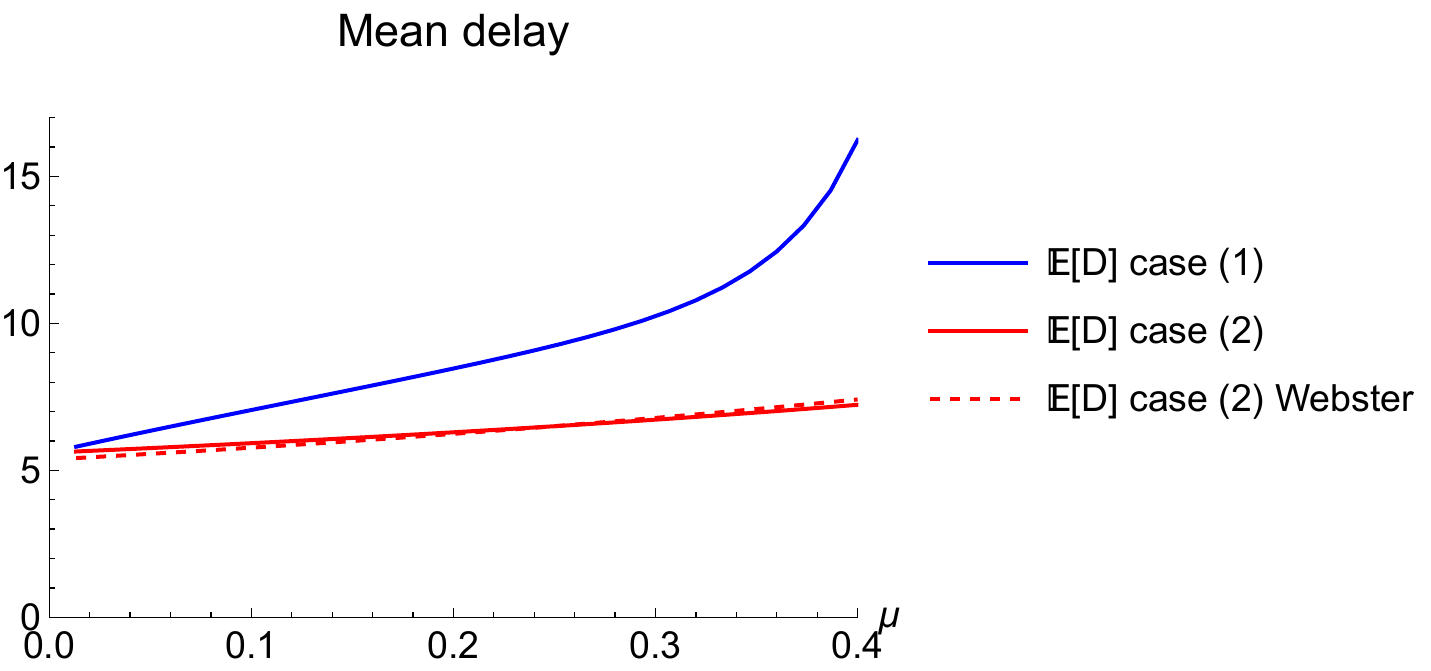}\\
			(b) \\[1ex]
		\caption{The total Poisson arrival rate, $\mu$, on the horizontal axis and in (a) the mean queue length at the end of an arbitrary time slot for the various cases and lanes where $\mathbb{E}[X^{t}]=\mathbb{E}[X^{(1)}]+\mathbb{E}[X^{(2)}]$ for case (2), and in (b) the mean delay for the various cases.}
		\label{f:single}
	\end{figure}

	We start with comparing the traffic performance of a single shared right-turn lane as in Figure~\ref{fig:vis2}(a), case ($1$), and a two-lane scenario where the turning vehicles and the straight-going vehicles are split as in Figure~\ref{fig:vis2}(b), case ($2$). We refer in the latter case to the lane which has right-turning vehicles as lane $1$ and to the other lane we refer as lane $2$. We assume that the arrival process is Poisson and that the arrival rate of turning vehicles, $\mu_1$, and straight-going vehicles, $\mu_2$, are the same in both cases. The total arrival rate of vehicles is $\mu=\mu_1+\mu_2$ in case~($1$). We choose $p_i=0.3$ for the shared right-turn lane, whereas in the two-lane case we have $p_i=1$ for lane $1$ and $p_i=0$ for lane $2$ and arrival rates $\mu_1=0.3\mu$ at lane $1$ and $\mu_2=0.7\mu$ at lane $2$. Further, we choose $q_i=1$, $g_1=8$, $g_2=20$, and $r=20$. We compute the mean queue length at the end of an arbitrary time slot for both lanes in case ($2$), denoted with  $\mathbb{E}[X^{(i)}]$ for lane $i$, and the total mean queue length at the end of an arbitrary time slot, denoted with $\mathbb{E}[X^{t}]$, and which equals $\mathbb{E}[X^{(1)}]+\mathbb{E}[X^{(2)}]$. For case ($1$) we denote the mean queue length at the end of an arbitrary time slot with $\mathbb{E}[X^t]$. The delay of an arbitrary car is denoted with $\mathbb{E}[D]$ for both cases ($1$) and ($2$). We study an example with various values of $\mu$ in Figure~\ref{f:single}.
	
%	\begin{table}[h!]
%		\begin{center}
%			\caption{The total Poisson arrival rate, $\mu$, on the horizontal axis and in (a) the mean queue length at the end of an arbitrary time slot for the various cases and lanes, and in (b) the mean delay for the various cases.}
%			\label{t:single}
%			\begin{tabular}{|c|cc|cHcHcc|}
%				\hline
%				& \multicolumn{2}{|c|}{Case ($1$)} & \multicolumn{6}{|c|}{Case ($2$)} \\ \hline
%				$\mu$ & $\mathbb{E}[X^t]$ & $\mathbb{E}[D]$ & $\mathbb{E}[X^{(1)}]$ & $\mathbb{E}[D^{(1)}]$ & $\mathbb{E}[X^{(2)}]$ & $\mathbb{E}[D^{(2)}]$ & $\mathbb{E}[X^t]$ & $\mathbb{E}[D]$ \\\hline
%				%$0.0400$ & $0.248$ & $6.20$ & $0.103$ & $8.56$ & $0.126$ & $4.51$ & $0.229$ & $5.72$ \\
%				$0.08$ & $0.542$ & $6.771$ & $0.208$ & $8.67$ & $0.260$ & $4.65$ & $0.468$ & $5.855$ \\
%				%$0.120$ & $0.880$ & $7.33$ & $0.316$ & $8.79$ & $0.403$ & $4.80$ & $0.719$ & $5.99$ \\
%				$0.16$ & $1.262$ & $7.889$ & $0.427$ & $8.90$ & $0.555$ & $4.96$ & $0.982$ & $6.140$ \\
%				%$0.200$ & $1.69$ & $8.46$ & $0.541$ & $9.02$ & $0.718$ & $5.13$ & $1.26$ & $6.29$ \\
%				$0.24$ & $2.179$ & $9.080$ & $0.658$ & $9.14$ & $0.892$ & $5.31$ & $1.550$ & $6.458$ \\
%				%$0.280$ & $2.74$ & $9.80$ & $0.778$ & $9.26$ & $1.08$ & $5.50$ & $1.86$ & $6.63$ \\
%				$0.32$ & $3.451$ & $10.79$ & $0.902$ & $9.39$ & $1.280$ & $5.72$ & $2.182$ & $6.818$ \\
%				%$0.360$ & $4.48$ & $12.4$ & $1.03$ & $9.52$ & $1.50$ & $5.94$ & $2.53$ & $7.02$ \\
%				$0.40$ & $6.496$ & $16.23$ & $1.159$ & $9.66$ & $1.733$ & $6.19$ & $2.892$ & $7.230$ \\\hline
%			\end{tabular}
%		\end{center}
%	\end{table}
	
	In Figure~\ref{f:single}, we can clearly see that the total mean queue length at the two lanes in case ($2$) is lower than the mean queue length at the single lane in case ($1$). This makes sense from various points of view: in case ($2$), we have twice as many lanes as in case ($1$), so we would expect a smaller total mean queue length in case ($2$). Moreover, in case ($1$), it might happen that straight-going vehicles are blocked. Such blockages cannot occur in case ($2$), as all turning traffic is on lane $1$ and all vehicles that go straight are on lane $2$. These two reasons are the main drivers for the performance difference in cases ($1$) and ($2$). From the point of view of the traffic performance, \textbf{it thus makes sense to split the traffic on a shared right-turn lane into two separate streams of vehicles on two lanes while assuming one lane available for departures in case ($1$) and two lanes in case ($2$)}. We observe similar results when looking at the mean delay and comparing cases ($1$) and ($2$).
	
	\begin{remark}
	We emphasized before that the blocking mechanism makes it impossible to use existing methods to analyze the queue lengths and delays. However, in this particular example we have chosen the parameter settings in such a way that case (2) \emph{can} be analyzed using existing methods. The reason is that we have two separate lanes, each with its own ``extreme'' blocking mechanism: lane 1 contains \emph{only} turning vehicles and \emph{all} of them are blocked during $g_1$. Essentially, this turns this lane into a regular FCTL queue with an extra long red period ($r + g_1$) and a shorter green period ($g_2$). Lane 2 contains only vehicles going straight, none of which are blocked. This means that this lane is essentially a regular FCTL queue as well. As a consequence, these two lanes can be analyzed separately using standard FCTL methods. When applying the method described in \cite{van2006delay}, the mean delay would be exactly the same as computed in Figure~\ref{f:single}(b). Moreover, this means that we can also use Webster's well-known approximation for the mean delay for case (2). This has also been visualized in Figure~\ref{f:single}(b) and, indeed, the approximation is remarkably accurate. Still, we stress that this is only possible because we have chosen an extreme blocking mechanism ($q_i=1$) in combination with Poisson arrivals (Webster's approximation only works for Poisson arrival processes).

	\end{remark}
	
	\subsubsection{Two lanes for the shared right-turn}
	
	Now we turn to an example where we still have two dedicated lanes as in case ($2$) of the previous example, one for turning traffic and one for straight-going traffic, see Figure~\ref{fig:vis2}(b), but we compare it with a two-lane example where the vehicles mix, see Figure~\ref{fig:vis2}(c). All turning vehicles will be on lane $1$, but we allow some straight-going traffic to be present on lane $1$ too. Lane $1$ is thus a shared right-turn lane. On lane $2$, we only have vehicles that are heading straight. This could, e.g., model a scenario in which some straight-going vehicles desire to take a specific lane,  strategically anticipating on an upcoming exit. Anticipation in lane changing behaviour is more generally investigated in e.g.~\cite{choudhury2013modelling} in urban scenarios. We could adapt the value of $p$ depending on this number of strategic vehicles. In order to make a comparison between the various cases that we study and that is as fair as possible, we assume the following: the total arrival rate and the fraction of turning vehicles are the same.
	
\begin{figure}[ht!]
		\centering
			\includegraphics[width=0.65\textwidth]{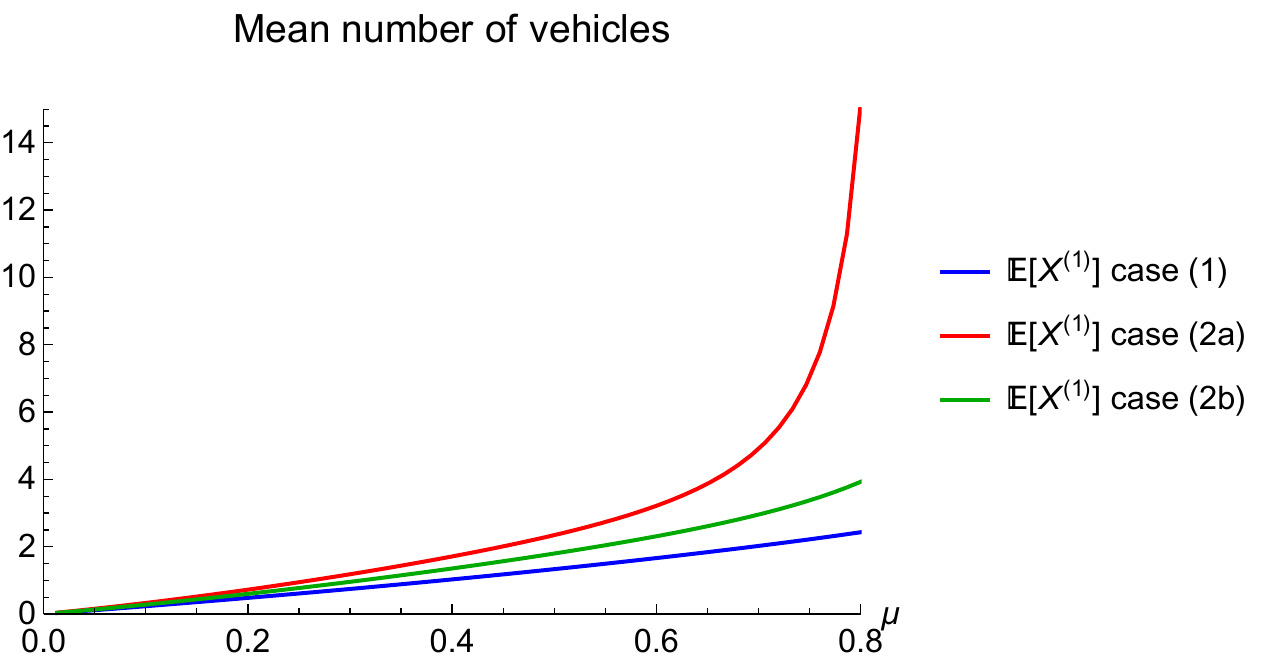}\\
			(a)\\[1ex]
			\includegraphics[width=0.65\textwidth]{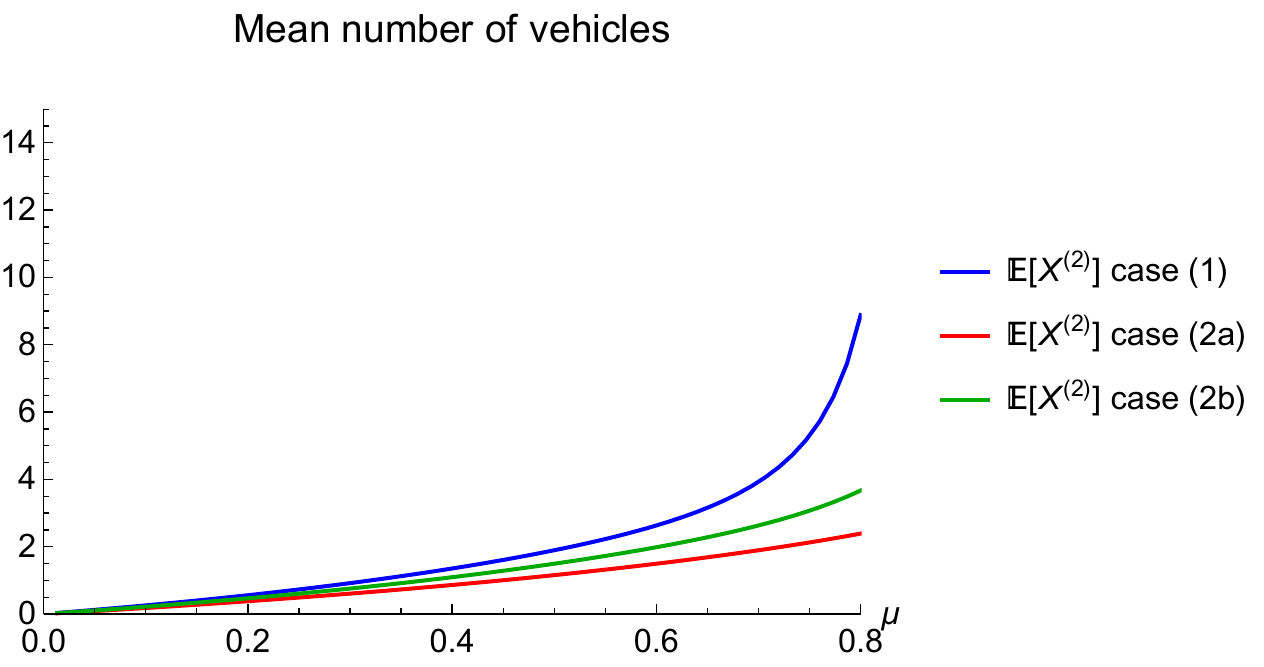}\\
			(b)\\[1ex]
			\includegraphics[width=0.65\textwidth]{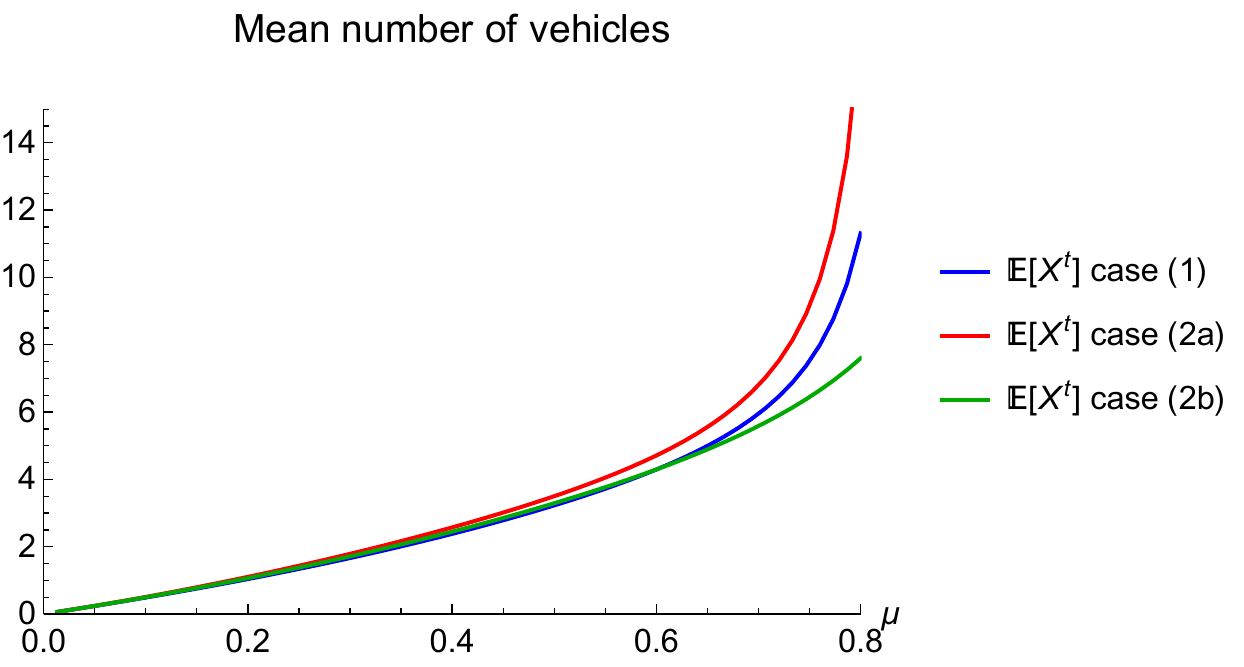}\\
			(c) \\[1ex]
		\caption{The total Poisson arrival rate, $\mu$ and the mean queue length at the end of an arbitrary time slot for the various cases, split among lane $1$ (a), lane $2$ (b) and the total among the two lanes (c) where $\mathbb{E}[X^{t}]=\mathbb{E}[X^{(1)}]+\mathbb{E}[X^{(2)}]$.}
		\label{f:split}
	\end{figure}

	We assume that the probability that an arbitrary vehicle is a turning vehicle is $0.3$ and we vary the total Poisson arrival rate $\mu$ to study the influence of the strict splitting of the turning vehicles. In case ($1$), we thus have an arrival rate at the right-turning lane that satisfies $\mu_1=0.3\mu$, whereas on the other lane we have an arrival rate $\mu_2 = 0.7\mu$. At lane $1$ we have $p_i=1$ and at lane $2$ we have $p_i=0$. In case ($2$) we distinguish between two subcases. In subcase ($2$a) we assume that the total arrival rate at both lanes is the same and thus $\mu_1=\mu_2=0.5\mu$. In subcase ($2$b), we assume that the arrival rate is split in the ratio $2:3$, so $\mu_1= 0.4\mu$ and $\mu_2=0.6\mu$. This implies that in subcase ($2$a) we choose $p_i=0.6$ (the fraction of turning vehicles is then $p\mu_1=0.6\cdot 0.5\mu=0.3\mu$) and in subcase ($2$b) we choose $p_i=0.75$ (the fraction of turning vehicles is then $p\mu_1=0.75\cdot0.4\mu=0.3\mu$), to make sure that we match the number of turning vehicles in case ($1$). Further, we choose $q_i=1$, $g_1=8$, $g_2=16$ and $r=16$. Then, we study the  mean queue length at the end of an arbitrary time slot of both lanes, $\mathbb{E}[X^{(1)}]$ and $\mathbb{E}[X^{(2)}]$, and the total average mean queue length at  the end of an arbitrary time slot, denoted with $\mathbb{E}[X^t]$. We obtain Figure~\ref{f:split}.
	
	In Figure~\ref{f:split}, we see only small differences in the total mean queue lengths at  the end of an arbitrary time slot for low arrival rates. At both lanes, there are few vehicles in the queue. This is different for the examples in Figure~\ref{f:split} with a higher arrival rate. In all examples for case ($1$) we see that the mean queue length at lane $2$, the straight-going traffic lane, is higher than for lane $1$. This is due to the relatively high fraction of vehicles that \emph{have} to use lane $2$ due to the strict splitting between turning and straight-going vehicles. In some sense, lane $1$, which only has turning vehicles, has overcapacity that cannot be used for the busier lane $2$ with only straight-going traffic. This is different for the other two cases, where the traffic is split more evenly across the two lanes. As one would expect, the longest queue in subcase ($2$a) is present at lane $1$, as the arrival rate at both lanes is the same and because vehicles are only blocked at lane $1$, the shared right-turn lane. This points towards another potential improvement and this is found in subcase ($2$b) where we balance the arrival rate differently. \emph{The right balance leads to a more economic use of both lanes and, hence, also the best performance} in this example when looking at $\mathbb{E}[X^t]$.
	
	The results in Figures~\ref{f:single} and \ref{f:split} might seem conflicting at a first glance, but they are not. In the case of a single, shared right-turn lane as in Figure~\ref{f:single}, we see a higher mean queue length than for the two dedicated lanes case in Figure~\ref{f:single}. This is the other way around in Figure~\ref{f:split} (considering case ($2$b)). This is mainly explained by the fact that in case ($2$b) in Figure~\ref{f:split}, we have two lanes and thus twice as many potential departures as in case ($1$) in Figure~\ref{f:single}. This is one of the main factors in the explanation of the differences in the mean performance between the examples studied in Figures~\ref{f:single} and \ref{f:split}.
	
	The two examples in this subsection tell us that \textbf{a separate or dedicated lane for turning traffic does not necessarily improve the traffic flow}. The intuition behind this is that a dedicated lane might have overcapacity which is not employed (e.g. in the case of an asymmetric load on both lanes). This issue is less present when the two dedicated lanes are turned into two lanes, one exclusively for straight-going traffic and one shared lane. This is confirmed by our simulations. As such, \textbf{an in-depth study is needed to obtain the best layout of the intersection and the best traffic-light control}. As a side-remark, we surpass the possibility here that in Figure~\ref{f:split}, case ($1$), we might control the two lanes in a different way, e.g. by prolonging the green period for one of the lanes. This is not possible in cases ($2$a) and ($2$b).

	\section{Conclusion and discussion}\label{sec:discussion}
	
	In this paper, we have established a recursion for the PGFs of the queue-length distribution at the end of each slot which can be used to provide a full queue-length analysis of the bFCTL queue with multiple lanes. This is an  extension of the regular FCTL queue so that we can account for temporal blockages of vehicles receiving a green light, for example because of a crossing pedestrian at the turning lane or because of a (separate) bike lane, and to account for a vehicle stream that is spread over multiple lanes. These features might impact
	the traffic-light performance as we have shown by means of various numerical examples. The blocking of turning vehicles and the number of lanes corresponding to a vehicle stream therefore has to be taken into account when choosing the settings for a traffic light.
	
	We briefly touched upon how one should design the layout of an intersection. Interestingly, it might be suboptimal to have a dedicated lane for turning traffic. It seems that mixing turning and straight-going traffic has benefits over a strict separation of those two traffic streams when there are two lanes for this turning and straight-going traffic. We advocate a further investigation into the influence of separating or mixing different streams of vehicles in front of traffic lights. It might be possible to find the optimal division of straight-going and turning vehicles over the various lanes, e.g. by enumerating several possibilities. A more structured optimization seems difficult because of the intricate expressions involved, but would definitely be worthwhile to investigate. Some research on the splitting of different traffic streams has already been done in e.g.~\cite{kikuchi2007lengths,tian2006probabilistic,wu1999capacity} and \cite{zhang2008modeling} and the present study can be seen as an alternative way of modelling the situation at hand.

	A possible extension of the results on the bFCTL queue is a study of (the PGF of) the delay distribution. We have refrained from deriving the delay distribution because of its (notational) complexity. Using proper conditioning, one can obtain (the PGF of) the delay distribution for the bFCTL queue.
	
	The work in~\cite{shaoluen2020random}, in which a simulation study of a similar model is performed, has been a source of inspiration for the study in this paper. %With the established model, there is no need to resort to simulations anymore.
	There are some extensions possible when comparing our work with~\cite{shaoluen2020random}. We e.g. did not study the influence of start-up delays as is done in~\cite{shaoluen2020random}. Investigating such start-up delays at the beginning of the green period is easily done in our framework: we simply need to adjust the $Y_i$ for the first few slots. Another approach to deal with start-up delays is presented in~\cite{maesnetworks}. Start-up delays which depend on the blocking of vehicles and different slot lengths for different combinations of turning/straight going vehicles, are harder to tackle. One could e.g. introduce additional states (besides states $u$ and $b$) to deal with this. Although the developed recursion does not directly allow for such a generalization, it seems possible to account for this at the expense of a more complex recursion. For the ease of exposition, we have refrained from doing so and we leave a full study on this topic for future research. %At the same time, we are already able to rigorously confirm some of the conclusions drawn by~\cite{shaoluen2020random}.
	
	%The bFCTL queue calls for further generalizations besides the ones mentioned in the previous paragraph. For example, instead of two full lanes, e.g. one for straight-going and one for turning traffic, we could also consider a single lane which splits into two lanes just before the intersection in such a way that some, say $N$, turning vehicles may accumulate on a separate lane. Such a small separate lane is often referred to as a turning bay. The $N$ vehicles on the turning bay would not block straight-going traffic in any way (because they are on a separate lane), but if there would be $N$ vehicles at the turning bay and another turning vehicle arrives, also the vehicles heading straight will be blocked. It would be interesting to study such a model and gain insight into the benefits of such a turning bay. %It seems that such a model is tractable and can be analyzed using similar techniques as the ones presented in this chapter.
	
	A further possible extension of the bFCTL queue would be to consider different blocking behaviours: instead of e.g. a fixed probability $q_i$ for each slot $i$, a more general blocking process might be considered. For example, if there are no pedestrians during slot $i$ for the model depicted in Figure~\ref{fig:vis}(b), then the probability that there are also no pedestrians in slot $i+1$, might be relatively high. In other words, there might be \emph{dependence} between the various slots when considering the presence of pedestrians. We gave an example where there is dependence between the current and the next slot, but it is also possible to consider such dependencies among more than two slots. It is worthwhile to investigate generalizations of the blocking process in order to further increase the general applicability of the bFCTL queue with multiple lanes.
	
	Another generalization for the blocking mechanism, is to block only a \emph{part} of the $m$ vehicles that are at the head of the queue. Indeed, we restrict ourselves to the cases where either all vehicles in a batch of size $m$ are blocked (or not). In various real-life examples, it might be the case that only part of the $m$ vehicles are blocked. It would be interesting to investigate whether such a model can be analyzed. Further, a situation with ``a right turn is always permitted'' scenario might be investigated. In such a case, right-turning vehicles are always free to turn, but might be blocked by straight-going vehicles in front them, which have to wait for a red traffic light, or are blocked by pedestrians. Straight-going vehicles might be blocked by turning traffic waiting for pedestrians. It seems that such a case, at the expense of additional complexity, can be tackled by a similar type of recursion as the one that is developed in this paper by extending and generalizing the blocking mechanism (and, thus, the recursion) to the red period.
	
	%Another topic for future research is to modify the bFCTL queue in such a way that it enables a \emph{joint} analysis of two lanes with one dedicated lane for vehicles heading straight and one shared lane with both turning and straight-going vehicles. Such a case is not covered by the bFCTL with multiple lanes, as it seems that in this extension one needs to take into account how many vehicles of both types (i.e. turning and straight-going vehicles) there are. This might lead to a two-dimensional queueing model rather than the one-dimensional one considered in the current paper.
	
	\paragraph{Discussion.} We end this paper with a discussion on its practical applicability. Although we have extended the standard model for traffic signals with fixed settings, there are still quite some possible improvements, as discussed in the above paragraphs. Still, to the best of our knowledge, this paper is the first to present analytical results for traffic intersections with blocking mechanisms, based on a queueing theoretic approach. Note that standard formulas like Webster's approximation for the mean delay \cite{webster1958traffic} cannot be used in these situations. From a practical point of view, the most relevant extension to the current analysis would be to deal with start-up delays that depend on the blocking of vehicles. One way to do this, is by considering different slot lengths for different combinations of turning/straight going vehicles, inspired by an analysis in ~\cite{maesnetworks}. This would make it possible to compute a saturation flow adjustment factor due to the right-turning movements at shared lane conditions (see also Biswas et al. \cite{biswas2018}).
	
	Finally, we also advocate an investigation whether the bFCTL queue with a vehicle-actuated mechanism (rather than the fixed green and red times that we consider) results in a tractable model. 
	
	\paragraph*{Acknowledgements} We would like to thank Onno Boxma for several interesting discussions that a.o. improved the readability of this manuscript. We are also thankful to Joris Walraevens who suggested the current exposition of the PGF recursion and to the reviewers who suggested several improvements of the paper.
	
	\paragraph*{Funding} The work in this paper is supported by the Netherlands Organization for Scientific Research (NWO) under grant number 438-13-206.
	
	\paragraph*{Disclosure statement} The authors report there are no competing interests to declare.

	\bibliographystyle{tfcad}
	\bibliography{yourbib} % if more than one, comma separated
	
\appendix
\newpage
\section*{Appendix}
\section{Simulation of the bFCTL queue}\label{a:simulation}
%\section{}

All numerical results in this paper have been obtained by implementing the bFCTL analysis in Mathematica. To validate these results, we have written a simulation of the model in Python. This simulation is a slightly more general version of the Python simulation used in the paper by Huang et al. \cite{shaoluen2020random}. In this paper, we wanted to show that their model can (also) be analysed using methods from queueing theory. Although the main message of the current paper is that this exact analysis is the preferred solution method, we want to give more insight in the simulation used to validate our results. For this reason, we include the most relevant parts of the Python code to simulate the bFCTL queue (see Listing~\ref{lst:pythonsim}).

\begin{lstlisting}[caption={Simulation of the bFCTL queue in Python.},
float=htb,label=lst:pythonsim,gobble=4]
    arrs = random.poisson(arrRate, (c, ncycles))    # random arrivals
    noBlocks = zeros((c, ncycles))
    noDeps = zeros((c, ncycles))
    
    for i in range(ncycles):
        for j in range(c):
            noBlocks[j,i] = random.binomial(1, q[j]) # blockages
            noDeps[j,i] = random.binomial(1, p[j])   # turning cars
    
    X = 0         # current queue length
    Xg = []       # queue length at end of green period
    Xi = zeros(c) # queue length at end of each time slot
    slot = 0      # slot counter

    blocked = False
    for i in range(ncycles):
        for j in range(c):
            # handle arrivals and departures if blocked
            if slot < g1:   # g1 period
                if blocked:   # blocked during g1 period
                    if noBlocks[j, i] == 0:  # with probability 1-q[j] 
                        blocked = False      # blockage resolved
                    else:
                        X += arrs[j, i]  # queue remains blocked
                else:       # not blocked during g1 period
                    if noDeps[j, i] == 1:   # with probability p[j]
                        if noBlocks[j, i] == 1:   # with prob. q[j]
                            blocked = True  # block turning vehicles
                            X += arrs[j, i] #  arrivals join the queue
            elif slot < g1 + g2:    # g2 period
                blocked = False     # blockages are always resolved in this period
            else:
                X += arrs[j, i]     # red period
    
            # handle arrivals and departures if NOT blocked (or blockage was resolved)
            if slot < g1 + g2 and not blocked:  
                if X < m:   # all vehicles can depart
                    X = 0   # and no new ones will arrive
                else:
                    X += arrs[j, i] - m  # m vehicles depart, new ones arrive
                    if X < 0:            # in this case, the queue becomes empty
                        X = 0
                
            if slot == g1+g2-1:     # at end of green period, store the queue length X_g
                Xg.append(X)
            Xi[slot] += X
            time += 1
            slot += 1
            if slot == c:           # reset slot number at end of cycle
                slot = 0
                
    Xi = Xi / ncycles
    print(Xi)         # Print the mean queue length at end of time slot
    print(mean(Xg))   # Print the mean queue length at end of green period
\end{lstlisting}

To show the accuracy of the simulation, we have repeated the experiment of Subsection~\ref{sec:influenceofturning}. In more detail, we have run the bFCTL simulation for the example shown in Figure~\ref{fig:ex3xk}(b), with Poisson arrivals with rate 0.39, $g_1=2$, $g_2=4$, $r=4$, $q_i=1$ and $m=1$. The exact results (using our theoretical analysis) and the simulation results are given in Table~\ref{tab:simresults}. The confidence intervals are based on 100 runs of 10,000 cycles each. Indeed, the simulation is accurate and confirms the correctness of the formulas derived in this paper. However, to obtain this level of accuracy, the simulations ran for almost two minutes, whereas the analytical results are obtained in just a few seconds. Admittedly, surely the efficiency of the code in Listing~\ref{lst:pythonsim} can be improved and more time would be gained by running the simulations in parallel on multiple cores. Still, the analytical methods will always outperform the simulation in terms of computation time and accuracy.

\begin{table}[ht]
    \centering
    \begin{tabular}{c|ccc|ccc}
    \multicolumn{7}{c}{$\mathbb{E}[X_i]$}\\
    \hline
 & \multicolumn{3}{c|}{$p_i=0$} & \multicolumn{3}{c}{$p_i = 0.6$} \\
$i$ & Exact & Sim CI Lower & Sim CI Upper & Exact & Sim CI Lower & Sim CI Upper \\
\hline
1 &1.297   & 1.289  & 1.298 & 3.901  &3.852     &  3.927\\
2 &0.926 & 0.917  & 0.925 & 4.148  &  4.102   &    4.178 \\
3 &0.657 & 0.649  & 0.657 & 3.610  &  3.562   &    3.638 \\
4 &0.465 & 0.458  & 0.465 & 3.126  &  3.078   &    3.154 \\
5 &0.329 & 0.323  & 0.330 & 2.699  &  2.652   &    2.726 \\
6 &0.233 & 0.228  & 0.233 & 2.325  &  2.280   &    2.352 \\
7 &0.623 & 0.617  & 0.623 & 2.715  &  2.669   &    2.742 \\
8 &1.013  & 1.007  & 1.014 & 3.105  & 3.059    &   3.133 \\
9 &1.404  & 1.396  & 1.405 & 3.495  & 3.448    &   3.522 \\
10& 1.793 & 1.785 &  1.794& 3.885  &  3.838   &    3.912 \\
    \end{tabular}
    \caption{Simulation results for the bFCTL queue with Poisson arrivals with rate 0.39, $g_1=2$, $g_2=4$, $r=4$, $q_i=1$ and $m=1$. The confidence intervals are based on 100 runs of 10,000 cycles each. }
    \label{tab:simresults}
\end{table}

\end{document}